\input amstex
\documentstyle{amsppt}
\hsize=16cm
\vsize=23cm

\topmatter
\title
Generalized Jacobian Rings for Open Complete Intersections
\endtitle
\author{Masanori Asakura and Shuji Saito}
\endauthor

\address
{Graduate School of Mathematics,
Kyushu University 33
FUKUOKA 812-8581, JAPAN}
\endaddress
\address
{e-mail: asakura\@math.kyushu-u.ac.jp}
\endaddress
\address
{Graduate School of Mathematics,
Nagoya University, Chikusa-ku, NAGOYA, 464-8602, JAPAN}
\endaddress
\address
{e-mail: sshuji\@msb.biglobe.ne.jp}
\endaddress

\endtopmatter
\document
\NoBlackBoxes

\input amstex
\documentstyle{amsppt}
\hsize=16cm
\vsize=23cm

\centerline{\bf Abstract}
\medbreak

  In this paper, we develop the theory of Jacobian rings of open complete
intersections, which mean a pair $(X,Z)$ where $X$ is a smooth complete
intersection in the projective space and and $Z$ is a simple normal crossing
divisor in $X$ whose irreducible components are smooth hypersurface sections
on $X$. Our Jacobian rings give an algebraic description of the cohomology of
the open complement $X-Z$ and it is a natural generalization of the Poincar\'e
residue representation of the cohomology of a hypersurface originally invented
by Griffiths. The main results generalize the Macaulay's duality theorem and
the Donagi's symmetrizer lemma for usual Jacobian rings for hypersurfaces.
A feature that distinguishes our generalized Jacobian rings from usual ones is
that there are instances where duality fails to be perfect while the defect can
be controlled explicitly by using the defining equations of $Z$ in $X$.
Two applications of the main results are given: One is the infinitesimal
Torelli problem for open complete intersections. Another is an explicit bound
for Nori's connectivity in case of complete intersections.
The results have been applied also to study of algebraic cycles in several
other works.

\vskip 20pt

\input amstex
\documentstyle{amsppt}
\hsize=16cm
\vsize=23cm
\document

\head \bf Contents
\endhead

\vskip 10pt
\roster
\item"\S0"
Introduction
\item"\S1"
Jacobian rings for open complete intersections
\item"\S2"
Green's Jacobian rings
\item"\S3"
Proof of Theorem(I)
\item"\S4"
A vanishing lemma
\item"\S5"
Proof of Theorem(II)
\item"\S6"
Proof of Theorem(II')
\item"\S7"
Proof of Theorem(III)
\item"\S8"
Infinitesimal Torelli for open complete intersections
\item"\S9"
Explicit bound for Nori's connectivity
\item"\hbox{ }"
References

\endroster

\vskip 20pt

\input amstex
\documentstyle{amsppt}
\hsize=16cm
\vsize=23cm

\def\Th#1.{\vskip 6pt \medbreak\noindent{\bf Theorem(#1).}}
\def\Cor#1.{\vskip 6pt \medbreak\noindent{\bf Cororally(#1).}}
\def\Pr#1.{\vskip 6pt \medbreak\noindent{\bf Proposition(#1).}}
\def\Lem#1.{\vskip 6pt \medbreak\noindent{\bf Lemma(#1).}}
\def\Rem#1.{\vskip 6pt \medbreak\noindent{\it Remark(#1).}}
\def\Fact#1.{\vskip 6pt \medbreak\noindent{\it Fact(#1).}}
\def\Claim#1.{\vskip 6pt \medbreak\noindent{\it Claim(#1).}}
\def\Def#1.{\vskip 6pt \medbreak\noindent{\bf Definition\bf(#1)\rm.}}

\def\qwith{\quad\hbox{with }}
\def\mathrm#1{\rm#1}

\def\isom{@>\cong>>}

\def\Ext{{\text{\rm{Ext}}}}

\def\dim{{\operatorname{dim}}}

\def\Coker{{\text{\rm Coker}}}
\def\dim{\hbox{\rm dim}}
\def\det{\hbox{\rm det}}

\def\Im{\hbox{\rm Im}}
\def\Ker{\hbox{\rm Ker}}
\def\Coker{\hbox{\rm Coker}}
\def\min{\hbox{\rm min}}

\def\Hom{\hbox{\rm{Hom}}}

\def\sign{\hbox{\mathrm{sign}}}
\def\Res#1{\hbox{\mathrm{Res}}_{#1}}

\def\P{{\Bbb{P}}}
\def\bP{{\Bbb{P}}}

\def\Q{{\Bbb{Q}}}

\def\cHom{{\Cal{H}}om}
\def\cL{{\Cal{L}}}
\def\L{{\Cal{L}}}
\def\cE{{\Cal{E}}}
\def\E{{\Cal{E}}}
\def\cF{{\Cal{F}}}
\def\cO{{\Cal{O}}}
\def\cX{{\Cal{X}}}

\def\cZ{{\Cal{Z}}}

\def\cM{{\Cal{M}}}

\def\sE#1#2#3{E_{#1}^{#2,#3}}

\def\l{\ell}
\def\d{{\bold{d}}}
\def\e{{\bold{e}}}

\def\lra{\longrightarrow}
\def\ra{\rightarrow}

\def\hra{\hookrightarrow}

\def\ot{\otimes}
\def\op{\oplus}

\def\us#1#2{\underset{#1}\to{#2}}
\def\os#1#2{\overset{#1}\to{#2}}

\def\PE{\P(\cE)}
\def\qaq{\quad\hbox{ and }\quad}
\def\qfor{\quad\hbox{ for }}
\def\qif{\quad\hbox{ if }}
\def\ld{\lambda}
\def\spa{\hbox{ }}
\def\scs{\spa : \spa}
\def\SS{\Sigma^*}
\def\sp{\hbox{}}

\def\Sp{\Sigma'}
\def\Sb{\overline{\Sigma}}
\def\SSp{\Sigma^{'*}}

\def\taub{\overline{\tau}}
\def\Lb{\overline{\Cal L}}
\def\Pb{\overline{\Bbb P}}
\def\Rp#1#2{B'_{#1}(#2)}
\def\Rb#1#2{\overline{B}_{#1}(#2)}
\def\eb{\overline{\bold e}}
\def\pg#1{{\frak S}_{#1}}
\def\ul#1{\underline{#1}}

\def\ba{\bold{a}}
\def\bb{\bold{b}}
\def\ua{\underline{a}}
\def\ub{\underline{b}}

\def\ue{\underline{e}}
\def\ud{\underline{d}}

\def\md{\delta_{\text{min}}}
\def\Om{\Omega}
\def\ld{\lambda}
\def\Wb{\overline{W}}

\def\FS#1{F_S^{#1}}
\def\WS#1{\Omega_S^{#1}}

\def\otO{\otimes_{\cO}}
\def\cU{{\Cal U}}

\def\onab{\overline{\nabla}}

\def\HcU#1#2{\HO^{#1,#2}(\cU/S)}
\def\HccU#1#2{H_{\cO,c}^{#1,#2}(\cU/S)}

\def\HcU#1#2{\HO^{#1,#2}(\cU/S)}

\def\HX#1#2{H^{#1}(X,{#2})}

\def\Xx{X_x}
\def\Zx{Z_x}
\def\Ux{U_x}

\def\Zxx{Z_{x}}
\def\Zst{Z}
\def\cZst{\cZ}
\def\Zast{Z^{(\alpha)}}
\def\Wst{W}
\def\Yst{Y}

\def\rocXZ{\kappa_{(\cX,\cZ)}}

\def\roolog{\kappa_0^{log}}

\def\roxlog{\kappa_x^{log}}

\def\TS{\Theta_S}
\def\TxS{T_x S}
\def\ToS{T_0 S}

\def\WS#1{\Omega_S^{#1}}
\def\WX#1{\Omega_X^{#1}}

\def\TX{T_X}

\def\TXZx{T_{\Xx}(-\log \Zxx)}

\def\WXZ#1{\Omega_X^{#1}(\log \Zst)}
\def\WXZa#1{\Omega_X^{#1}(\log \Zast)}
\def\WcXZ#1{\Omega_{\cX/S}^{#1}(\log \cZst)}
\def\WcXkZ#1{\Omega_{\cX/k}^{#1}(\log \cZst)}
\def\WcXCZ#1{\Omega_{\cX/\Bbb C}^{#1}(\log \cZst)}
\def\WPnY#1{\Omega_{\Bbb P^n}^{#1}(\log \Yst)}
\def\TXZ{T_X(-\log \Zst)}
\def\TcXZS{T_{\cX/S}(-\log \cZst)}

\def\fXZ#1{\phi_{X,Z}^{#1}}
\def\cXZ{\psi_{X,Z}}

\def\cXZx{\psi_{\Xx,\Zx}}

\def\eXZ{\eta_{X,Z}}

\def\ccXZS{c_S(\cX,\cZ)}

\def\dpXZ#1{d\rho^{#1}_{X,Z}}
\def\scs{\hbox{ }:\hbox{ }}
\def\onab{\overline{\nabla}}
\def\tonab{@>\onab>>}

\def\spa{\hbox{ }}

\def\otk{\otimes_k}
\def\Pol{P}
\def\h#1#2{h_{#1}(#2)}

\def\Szar{S_{zar}}

\head \S0. Introduction. \endhead
\vskip 8pt

The purpose of this paper is to develope the theory of
\it Jacobian rings of open complete intersections. \rm
Here, by ``open complete intersection" we mean a pair
$(X,Z=\underset{1\leq j\leq s}\to{\cup} Z_j)$ where $X$ is a smooth complete
intersection in $\Bbb P^n$ and $Z_j\subset X$ is a smooth hypersurface section
such that $Z$ is a simple normal crossing divisor in $X$.
Our Jacobian rings give an algebraic description of the cohomology of
$X\setminus Z$ and it is a natural generalization of the Poincar\'e residue
representation of the cohomology of a hypersurface that played a significant
role in the work of Griffiths [Gri].
The fundamental results on the generalized Jacobian rings have been stated
without proof in [AS] where it is applied to the Beilinson's Hodge and Tate
conjecture for open complete intersections (see [A], [MSS] and [SaS] for other
applications of the generalized Jacobian rings).
The main results are stated in \S1. The proofs occupying \S2 through \S7
are based on the basic techniques in computations of Koszul cohomology
developed by M. Green ([G1] and [G2]).
Two applications of the results in \S1 are given in \S8 and \S9.
\vskip 4pt

In \S8 we study the infinitesimal Torelli problem for open complete
intersection $(X,Z)$ as an application of the duality theorem for the
generalized Jacobian rings. It concerns the injectivity of the following map
$$ \HX 1 {\TXZ} \to \underset{1\leq q\leq m}\to{\bigoplus}
\Hom(\HX {m-q} {\WXZ q},\HX {m-q+1} {\WXZ {q-1}}) $$
induced by the cup product and the contraction
$\WXZ q \ot \TXZ \to \WXZ {q-1}$.
Here $\WXZ q$ is the sheaf of differential $q$-forms on $X$
with logarithmic poles along $Z$ and $\TXZ$ is the $\cO_X$-dual of $\WXZ 1$.
We show that the map is injective under a mild numerical assumption.
Since the above map is interpreted as the derivative of the period map
from an appropriate moduli space of isomorphism classes of pairs $(X,Z)$
to the period domain (cf. [U2]),
it implies that the mixed Hodge structure on $H^m(X\setminus Z,\Bbb Q)$
determines $(X,Z)$ up to isomorphisms locally on the moduli space.
It is a generalization of the infinitesimal Torelli
for hypersurfaces due to Griffiths [Gri] and for complete intersections due
to Peters [P] and Usui [U1].
\vskip 4pt

In \S9 we prove the following result as an application of the symmetrizer lemma
for the generalized Jacobian rings.
We fix integers $r,s\geq 1$ and $d_1,\dots,d_r,e_1,\dots,e_s\geq 1$.
Let $S$ be a non-singular affine variety over $\Bbb C$ and
assume that we are given schemes over $S$
$$\Bbb P^n_S \hookleftarrow \cX \hookleftarrow
\cZ=\underset{1\leq j\leq s}\to{\cup} \cZ_j,\leqno(0-1)$$
whose fibers are open complete intersections.
Assume that the fibers of $\cX/S$ are smooth complete intersection of
multi-degree $(d_1,\dots,d_r)$ in $\Bbb P^n$ and that those of
$\cZ_j\subset\cX$ are smooth hypersurface section of degree $e_j$.
Write $\cU=\cX \setminus \cZ$.
We will introduce an invariant $\ccXZS$ that measures
the ``generality" of the family (0-1), or how many independent parameters $S$
contains.

\Th 0-1. \it Assuming $n-r\geq 2$, we have
$$ \align
&F^{t-n+r+1} H^t(\cU,\Bbb C)=0
\qif
s\leq n-r+2
\hbox{ and }
\md r\geq t+r+1+\ccXZS,\\
& F^{t-n+r+1} H^t(\cX,\cZ,\Bbb C)=0
\qif s=1
\hbox{ and }
\md r+e_1\geq t+r+1+\ccXZS.\\
\endalign$$
where $\md=\underset{1\leq i\leq r,1\leq j\leq s}\to{\min}
\{d_i,e_j\}$ and
$F^*$ denotes the Hodge filtration defined in [D1] and [D2].
\rm\vskip 6pt

The second vanishing of Th.(0-1) gives an explicit bound for Nori's
connectivity [N] in case of complete intersections in the projective space.
Nagel [Na2] has obtain similar degree bounds for complete intersections
in a general projective smooth variety.

\vskip 6pt

\vskip 20pt

\input amstex
\documentstyle{amsppt}
\hsize=16cm
\vsize=23cm

\head \S1. Jacobian rings for open complete intersections.
\endhead
\vskip 8pt

Throughout the whole paper, we fix integers $r,s \geq 0$ with $r+s\geq 1$,
$n\geq 2$ and $d_1, \cdots, d_r$, $e_1, \cdots, e_s \geq 1$. We put
$$
\d=\sum_{i=1}^{r}d_i, \quad \e=\sum_{j=1}^{s}e_j,\quad
\md=\underset{\underset{1\leq j\leq s}\to{1\leq i\leq r}}\to{\min}\{d_i,e_j\},
\quad d_{max}=\underset{1\leq i\leq r}\to{\max}\{d_i\},
\quad e_{max}=\underset{1\leq j\leq s}\to{\max}\{e_j\}.
$$
We also fix a field $k$ of characteristic zero.
Let $\Pol=k[X_0,\dots,X_n]$ be the polynomial ring over $k$ in $n+1$
variables.
Denote by $\Pol^l\subset \Pol$ the subspace of the homogeneous
polynomials of degree $l$.
Let $A$ be a polynomial ring over $\Pol$ with indeterminants
$\mu_1, \cdots, {\mu}_r$, $\lambda_1,\cdots, \lambda_s$.
We use the multi-index notation
$$\mu^{\underline{a}}=\mu_1^{a_1}\cdots\mu_r^{a_r} \qaq
 \lambda^{\underline{b}}=\lambda_1^{b_1}\cdots\lambda_s^{b_s} \qfor
\underline{a}=(a_1,\cdots,a_r) \in {\Bbb Z}^{\oplus r}_{\geq 0},\spa
\underline{b}=(b_1,\cdots,b_s) \in {\Bbb Z}^{\oplus s}_{\geq 0}.$$
For $q \in {\Bbb Z}$ and $\ell \in {\Bbb Z}$, we write
$$
 A_q(\ell)=
\us{\ba+\bb=q}{\oplus}
\Pol^{\ua\ud+\ub\ue+\ell} \cdot
\mu^{\underline{a}}\lambda^{\underline{b}}\quad
(\ba=\sum_{i=1}^r a_i,\sp
\bb=\sum_{j=1}^s b_j,\sp
\ua\ud=\sum_{i=1}^r a_id_i,\sp
\ub\ue=\sum_{j=1}^s b_je_j)
$$
By convention $A_q(\ell)=0$ if $q<0$.

\Def1-1.
For $\ul{F}=(F_1, \cdots, F_r)$, $\ul{G}=(G_1, \cdots, G_s)$ with
$F_i \in \Pol^{d_i}$, $G_j \in \Pol^{e_j}$, we define the {\bf Jacobian ideal}
$J(\ul{F},\ul{G})$ to be the ideal of $A$ generated by
$$\sum_{1\leq i\leq r}\frac{\partial F_i}{\partial X_k}\mu_i+
\sum_{1\leq j\leq s}\frac{\partial G_j}{\partial X_k}\lambda_j,
\quad F_i,\quad G_j\lambda_j \quad
(1 \leq i \leq r,\sp 1 \leq j \leq s,\sp 0 \leq k \leq n).$$
The quotient ring $B=B(\ul{F},\ul{G})=A/J(\ul{F},\ul{G})$ is called
the {\bf Jacobian ring of} ($\ul{F},\ul{G}$). We denote
$$
B_q(\ell) =B_q(\ell)(\ul{F},\ul{G})=A_q(\ell)/A_q(\ell) \cap J(\ul{F},\ul{G}).
$$
\vskip 6pt

\Def 1-2. Suppose $n\geq r+1$.
Let $\Bbb{P}^n=\text{\rm Proj }\Pol$ be the projective space over $k$.
Let $X\subset \P^n$ be defined by $F_1=\cdots=F_r=0$ and let
$Z_j\subset X$ be defined by $G_j=F_1=\cdots=F_r=0$ for $1 \leq j \leq s$.
We also call $B(\ul{F},\ul{G})$ the Jacobian ring of the pair
$(X,\Zst={\cup}_{1\leq j\leq s} Z_j)$ and denote
$B(\ul{F},\ul{G})=B(X,\Zst)$ and $J(\ul{F},\ul{G})=J(X,\Zst)$.

\vskip 6pt
In what follows we fix $\ul{F}$ and $\ul{G}$ as Def.(1-1) and assume
the condition
$$ \text{$F_i=0$ ($1\leq i\leq r$) and
$G_j=0$ ($1\leq j\leq s$) intersect transversally in $\Bbb{P}^n$.}
\leqno(1-1)$$
We mension three main theorems.
The first main theorem concerns with the geometric meaning of Jacobian rings.

\Th I. \it
Suppose $n\geq r+1$. Let $X$ and $Z$ be as Definition (1-2).
\vskip 4pt\noindent
(1)
For intergers $0 \leq q\leq n-r$ and $\ell\geq 0$ there is a natural
isomorphism
$$\fXZ q \scs B_q(\d+\e-n-1+\ell) \isom H^{q}(X, \WXZ {n-r-q}(\ell))_{prim}.$$
Here $\WXZ p$ is the sheaf of algebraic differential $q$-forms on $X$
with logarithmic poles along $\Zst$ and `$prim$' means the primitive part:
$$ H^{q}(X, \WXZ {p}(\ell))_{prim}=
\left\{\aligned
& \Coker(H^q(\bP^n,\Omega_{\bP^n}^q) \to H^{q}(X, \WX {q}))
\quad\text{ if $q=p$ and $s=\ell=0$,}\\
& H^{q}(X, \WXZ {p}(\ell))
\quad\text{ otherwise.}\\
\endaligned\right.$$
\vskip 4pt\noindent
(2)
There is a natural map
$$\cXZ : B_1(0) \longrightarrow H^1(X, \TXZ)_{alg}\subset H^1(X, \TXZ) $$
which is an isomorphism if $\dim(X)\geq 2$.
Here $\TXZ$ is the $\cO_X$-dual of $\WXZ 1$ and the group on the right hand
side is defined in Def.(1-3) below. The following map
$$
H^1(X, \TXZ) \otimes H^{q}(X, \WXZ p)
\longrightarrow H^{q+1}(X, \WXZ {p-1}).
$$
induced by the cup-product and the contraction
$\TXZ\ot\WXZ p\to\WXZ {p-1}$
is compatible through $\cXZ$ with the ring multiplication up to scalar.
\rm\vskip 6pt

Roughly speaking, the generalized Jacobian rings describe the infinitesimal
part of the Hodge structures of open
variety $X \setminus Z$,
and the cup-product with Kodaira-Spencer class
coincides with the ring multiplication up to non-zero scalar.
This result was originally invented by P. Griffiths in case of hypersurfaces
and generalized to complete intersections by Konno [K].
Our result is a further generalization.

\Def 1-3. \it Let the assumption be as in Th.(I).
We define $H^1(X, \TXZ)_{alg}$ to be
the kernel of the composite map
$$H^1(X,\TXZ)\to H^1(X,\TX) \to H^2(X,\cO_X),$$
where the second map is induced by the cup product with the class
$c_1(\cO_{X}(1))\in H^1(X,\WX 1)$ and the contraction
$\TX\otimes \WX 1 \to \cO_X$.
It can be seen that
$$ \dim_k(H^1(X, \TXZ)/H^1(X, \TXZ)_{alg})=\left\{
\aligned
1 & \text{ if $X$ is a $K3$ surface,}\\
0 & \text{ otherwise.}\\
\endaligned \right.
$$

\vskip 6pt\rm

The second main theorem is the duality theorem for the generalized
Jacobian rings.

\Th II. \it
(1)
There is a natural map (called the trace map)
$$\tau\sp:\sp B_{n-r}(2(\d-n-1)+\e) \to k.$$
Let
$$ \h p \ell\sp:\sp B_p(\d-n-1+\l) \ra B_{n-r-p}(\d+\e-n-1-\l)^*$$
be the map induced by the following pairing induced by
the multiplication
$$B_{p}(\d-n-1+\l) \otimes B_{n-r-p}(\d+\e-n-1-\l)\to  B_{n-r}(2(\d-n-1)+\e)
@>\tau>> k.$$
When $r>n$ we define $ \h p \ell$ to be the zero map by convention.
\vskip 4pt\noindent
(2)
The map $\h p \ell$ is an isomorphism in either of the following cases.
\roster
\item"$(i)$"  $s\geq 1$ and $p<n-r$ and $\ell<e_{max}$.
\item"$(ii)$"  $s\geq 1$ and $0\leq \ell\leq e_{max}$ and $r+s\leq n$.
\item"$(iii)$" $s=\ell=0$ and either $n-r\geq 1$ or $n-r=p=0$.
\endroster
\vskip 4pt\noindent
(3)
The map $\h {n-r} \ell$ is injective if $s\geq 1$ and $\ell<e_{max}$.
\vskip 6pt\rm

We have the following auxiliary result on the duality.

\Th II'. \it
Assume $n-r\geq 1$ and consider the composite map
$$ \eXZ \scs H^{0}(X, \WXZ {n-r})@>{(\fXZ 0)^{-1}}>>
B_0(\d+\e-n-1) @>{{\h {n-r} 0}^*}>>  B_{n-r}(\d-n-1)^*$$
where the second map is the dual of $\h {n-r} 0$.
Then $\eXZ$ is surjective and we have (cf. Def.(1-4) below)
$$ \Ker(\eXZ)=\wedge_X^{n-r}(G_1,\dots,G_s).$$
\rm\vskip 6pt

\Def 1-4. \it
Let $G_1,\dots,G_s$ be as in Def.(1-1) and let $Y_j\subset \P^n$ be
the smooth hypersurface defined by $G_j=0$.
Let $X\subset \P^n$ be a smooth projective variety such that
$Y_j$ ($1\leq j\leq s$) and $X$ intersect transversally.
Put $Z_j=X\cap Y_j$.
Take an integer $q$ with $0\leq q \leq s-1$.
For integers $1\leq j_1<\cdots< j_{q+1}\leq s$, let
$$ \omega_X(j_1,\dots,j_{q+1})\in H^{0}(X, \WXZ {q})
\quad (\Zst=\underset{1\leq j\leq s}\to{\Sigma} Z_j)$$
be the restriction of
$$ \sum_{\nu=1}^{q+1}(-1)^{\nu-1} e_{j_{\nu}}
\frac{dG_{j_1}}{G_{j_1}}\wedge\cdots\wedge\widehat
{\frac{dG_{j_\nu}}{G_{j_\nu}}}
\wedge\cdots\wedge{\frac{dG_{j_{q+1}}}{G_{j_{q+1}}}} \sp\in
H^0(\P^n,\WPnY q)$$
where
$\Yst=\underset{1\leq j\leq s}\to{\Sigma} Y_j\subset \P^n$.
We let
$$ \wedge_X^q(G_1,\dots,G_s)\subset H^{0}(X, \WXZ {q})$$
be the subspace generated by $\omega_X(j_1,\dots,j_{q+1})$.
For $1\leq j_1<\cdots< j_{q}\leq s-1$ we have
$$  e_s\cdot \omega_X(j_1,\dots,j_{q},s)=
\frac{dg_{j_1}}{g_{j_1}}\wedge\cdots\wedge
{\frac{dg_{j_{q}}}{g_{j_{q}}}}
\qwith g_j=(G_j^{e_s}/G_s^{e_j})_{|X}\in \Gamma(U,\cO_{U}^*)
\quad (U=X\setminus Z)$$
and $\omega_X(j_1,\dots,j_{q},s)$ with $1\leq j_1<\cdots< j_{q}\leq s-1$
form a basis of $\wedge_X^q(G_1,\dots,G_s)$.
\rm\vskip 6pt

Our last main theorem is the generalization of Donagi's symmetrizer lemma [Do]
(see also [DG], [Na1] and [N]) to the case of open complete intersections at
higher degrees.

\Th III. \it Assume $s\geq 1$.
Let $V \subset B_1(0)$ is a subspace of codimension $c\geq 0$.
Then the Koszul complex
$$
B_p(\l) \otimes \os{q+1}{\wedge}V \rightarrow
B_{p+1}(\l) \otimes \os{q}{\wedge}V \rightarrow
B_{p+2}(\l) \otimes \os{q-1}{\wedge}V
$$
is exact if one of the following conditions is  satisfied.
\roster
\item"$(i)$"
$p\geq 0$, $q=0$ and $\md p+\ell\geq c$.
\item"$(ii)$"
$p\geq 0$, $q=1$ and $\md p+\ell\geq 1+c$ and $\md(p+1)+\ell\geq d_{max}+c$.
\item"$(iii)$"
$p\geq 0$, $\md(r+p)+\ell\geq \d+q+c$, $\d+e_{max}-n-1>\ell\geq \d-n-1$
and either $r+s\leq n+2$ or $p\leq n-r-[q/2]$, where
$[*]$ denotes the Gaussian symbol.
\endroster
\rm\vskip 6pt

\Rem 1-1. \it
In case $q\geq 2$, $\l=\d-n-1$, $r+s>n+2$ and $p=n-r-1$,
the complex in Th.(III) is not injective but the cohomology is controlled
by motivic elements. We shall study it in a future paper.
\rm

\vskip 20pt

\input amstex
\documentstyle{amsppt}
\hsize=16cm
\vsize=23cm

\head \S2. Green's Jacobian rings \endhead
\vskip 8pt

Let the notation be as \S1.
In this section we study our Jacobian rings by using the method in
[G2] and [G3]. We write
$$\cE=\cE_0\op\cE_1 \qwith \cE_0=\os{r}{\us{i=1}{\op}}\cO(d_i)
\text{ and }
\cE_1=\os{s}{\us{j=1}{\op}}\cO(e_j)$$
where we denote $\cO=\cO_{\P^n}$.
We consider the projective space bundle
$$
\pi:\P:=\P(\cE)\lra \P^n.
$$
We let $\cL=\cO_{\P(\cE)}(1)$, the tautological line bundle.
We let
$$
\mu_i\in H^0(\P,\cL \otimes \pi^*\Cal O(-d_i))\qaq
\lambda_j\in H^0(\P,\cL \otimes \pi^*\Cal O(-e_j))$$
be the global section associated respectively to the effective divisors
$$\Bbb P (\us{\alpha\not=i}{\oplus}\cO(d_\alpha)\oplus \cE_1)
\hookrightarrow \PE
\qaq
\Bbb P(\cE_0 \oplus \us{\beta\not=j}{\oplus}\cO(e_\beta))
\hookrightarrow \PE.$$
Further we fix a global section
$$\sigma = \sum^r_{i=1} F_i \mu_i+\sum^s_{j=1} G_j \ld_j \in H^0(\P,\L),
\leqno (2-1)$$
and put
$$
\Q_i:\mu_i=0\subset \P, \quad \P_j:\ld_j=0  \subset \P, \quad
\Cal Z: \sigma=0 \subset \P,$$
$$
X_i:F_i=0 \subset \P^n, \quad Y_j:G_j=0 \subset \P^n.$$
We assume that
$$\bigcup_{1\leq i\leq r} X_i \cup \bigcup_{1\leq j\leq s}  Y_j\subset \P^n
\text{ is a simple normal crossing divisor}\leqno(2-2)$$
that implies that $\cZ$ is a nonsingular divisor in $\P$.
We will use the following divisors on $\P$
$$
\Q_*=\sum_{1\leq i\leq r}\Q_i
\qaq
\P_*=\sum_{1\leq j\leq s}\P_j$$
The following facts are well-known.

\Lem 2-1. \it Put $t=r+s$.
\roster
\item
We have the isomorphisms
$$
R^i\pi_*\Cal L^{\nu} \simeq
\left.\left\{\gathered
 S^{\nu}(\Cal E)\\
 \text{\det} \Cal E^*\ot S^{-\nu-t}(\Cal E^*)\\
  0
\endgathered\right.\qquad
\aligned
&\text{if $\nu\geq0$, $i=0$}\\
&\text{if $\nu\leq-t$, $i=t-1$}\\
&\text{otherwise}
\endaligned\right.
$$

\item
We have the isomorphisms
$$
H^q(\Bbb P,\Cal L^{\nu} \ot \pi^*\Cal V)\simeq
\left.\left\{\gathered
 H^q(\Bbb P^n, S^{\nu}(\Cal E) \ot \Cal V)\\
 H^{q-t+1}(\Bbb P^n, S^{-\nu-t}(\Cal E^*) \ot \text{\det} \Cal E^* \ot \Cal V)\\
 0
\endgathered\right.\qquad
\aligned
&\text{if $\nu \geq 0$}\\
&\text{if $\nu \leq -t$}\\
&\text{if $-t+1\leq \nu \leq -1$}
\endaligned\right.
$$
where $\Cal V$ is a vector bundle on $\Bbb P^n$.

\item
We have the commutative diagram with the exact horizontal sequences
(called the Euler sequences)
$$
\CD
0 @>>> \cO_{\P} @>>> \pi^*\cE^*\ot\cL @>>> T_{\P/\P^n} @>>> 0\\
@.@AA{=}A@AAA@AAA\\
0 @>>> \cO_{\P} @>>> \pi^*\cE_0^*\ot\cL\op \cO_{\P}^{\op s} @>>>
T_{\P/\P^n}(-\log\P_*) @>>>0\\
@.@AA{=}A@AAA@AAA\\
0 @>>> \cO_{\P} @>>> \cO_{\P}^{\op r+s} @>>>
T_{\P/\P^n}(-\log\P_*+\Q_*) @>>>0,\\
\endCD$$
where the middle vertical maps are given by the global sections
$$ \mu_i\in H^0(\P,\cL\ot\pi^*\cO(-d_i))
\qaq \lambda_j\in H^0(\P,\cL\ot\pi^*\cO(-e_j)).$$
\endroster
\rm\vskip 8pt

We introduce the sheaf of differential operators
of $\cL$ of order $\leq1$ as follows:
$$\align
\Sigma_{\cL}:={\Cal D}\text{\it iff}^{~\leq1}(\cL)
&=\{P\in {\Cal E}\text{\it nd}_{k}(\cL)~;~
Pf-fP\text{ is $\cO_{\P}$-linear }(f\in\cO_{\P})
\}\\
&\simeq \cL\ot D_{\P}^{\leq1}\ot \cL^*.
\endalign$$
By definition it admits the following exact sequence
$$
0\lra \cO_{\P}\lra \Sigma_{\cL}\lra T_{\P}\lra0
\tag{2-3}
$$
with the extension class
$$-c_1(\cL)\in \Ext^1(T_{\P},\cO_{\P})\simeq
\Ext^1(\cO_{\P},\Omega_{\P}^1\ot\cO_{\P})\simeq
H^1(\P,\Omega^1_{\P}).$$
Letting $U\subset \P^n$ be an affine subspace
and $x_1,\cdots,x_n$ be its coordinate,
$\varGamma(\pi^{-1}(U),\Sigma_{\cL})$ is
generated by the following sections
$$
\frac{\partial}{\partial x_i},~
\lambda_i\frac{\partial}{\partial \lambda_j},~
\lambda_i\frac{\partial}{\partial \mu_j},~
\mu_i\frac{\partial}{\partial \lambda_j},~
\mu_i\frac{\partial}{\partial \mu_j},~
\cO_{\P}\text{-linear maps}.
\tag{2-4}
$$
The section $\sigma$ defines a map
$$
j(\sigma):\Sigma_{\cL}\lra \cL,
\quad P\longmapsto P(\sigma),
$$
which is surjective by the assumption (2-2) and it gives rise to the
exact sequence
$$
0\lra T_{\P}(-\log \cZ)\lra \Sigma_{\cL}\os{j(\sigma)}{\lra} \cL\lra 0.
\tag{2-5}
$$

\Def 2-1. \it We define
$$ \Sigma_{\cL}(-\log \P_*)\subset \Sigma_{\cL} \qaq
\Sigma_{\cL}(-\log \P_*+\Q_*)\subset \Sigma_{\cL}$$
to be the inverse image of $T_{\P}(-\log \P_*)$ and $T_{\P}(-\log \P_*+\Q_*)$
respectively via the map in (2-3).
\rm\vskip 4pt

By the assumption (2-2) $j(\sigma)$ restricted on
$\Sigma_{\cL}(-\log \P_*+\Q_*)$ is surjective so that (2-5) gives rise to
the following commutative diagram of the exact sequences
$$\CD
0@>>>  T_{\P}(-\log \cZ) @>>> \Sigma_{\cL} @>{j(\sigma)}>>  \cL @>>> 0\\
@.@AAA@AAA@AA{=}A\\
0 @>>> T_{\P}(-\log \cZ+\P_*) @>>> \Sigma_{\cL}(-\log\P_*) @>{j(\sigma)}>>
 \cL @>>> 0\\
@.@AAA@AAA@AA{=}A\\
0 @>>> T_{\P}(-\log \cZ+\P_*+\Q_*) @>>>  \Sigma_{\cL}(-\log\P_*+\Q_*)
@>{j(\sigma)}>>   \cL @>>> 0.\\
\endCD
\tag{2-6}
$$
On the other hand (2-3) and the Euler sequences give rise to the
following commutative diagram of the exact sequences
$$
\CD
0 @>>> \pi^*\cE^*\ot\cL @>{\iota}>> \Sigma_{\cL}@>>> \pi^* T_{\P^n} @>>>0\\
@.@AAA@AAA@AA{=}A\\
0 @>>> \pi^*\cE_0^*\ot\cL\op \cO_{\P}^{\op s} @>>> \Sigma_{\cL}(-\log\P_*)
@>>> \pi^* T_{\P^n} @>>>0\\
@.@AAA@AAA@AA{=}A\\
0 @>>> \cO_{\P}^{\op r+s} @>>> \Sigma_{\cL}(-\log \P_*+\Q_*) @>>>
\pi^* T_{\P^n} @>>>0.
\endCD \tag{2-7}$$
Here $\iota$ is the sum of the maps
$\cL\ot\pi^*\cO(-d_i)\to\Sigma_{\cL}$ and
$\cL\ot\pi^*\cO(-e_j)\to\Sigma_{\cL}$
given by
$$ \frac{\partial}{\partial \mu_i}
\in H^0(\P,\cL^{-1}\ot\pi^*\cO(d_i)\ot\Sigma_{\cL})
\qaq
 \frac{\partial}{\partial \lambda_j}
\in H^0(\P,\cL^{-1}\ot\pi^*\cO(e_j)\ot\Sigma_{\cL})$$
respectively. The left vertical maps are given by the global sections
$$ \mu_i\in H^0(\P,\cL\ot\pi^*\cO(-d_i))
\qaq \lambda_j\in H^0(\P,\cL\ot\pi^*\cO(-e_j)).$$

The following lemma is straightforward from the definition.
\vskip 4pt

\Lem 2-2. \it For integers $\ell$ and $q$, we have
$$A_q(\l)=H^0(\cL^q\ot \pi^*\cO(\l))
\qaq B_q(\l)=A_q(\l)/J(\Sigma_{\cL}(-\log\P_*)).$$
where $J(\Sigma_{\cL}(-\log\P_*))\subset A_q(\l)$ is
the image of the map
$$
j(\sigma)\ot 1:H^0(\Sigma_{\cL}(-\log\P_*)\ot \cL^{q-1}\ot\pi^*\cO(\l))
\ra H^0(\cL^q\ot\pi^*\cO(\l)).
$$
\rm\vskip 6pt

We define the sheaf of differential operators
of $\cE$ of order $\leq1$ as follows:
$$\align
\Sigma'_{\cE}:={\Cal D}\text{\it{iff}}^{~\leq1}(\cE)
&=\{P\in {\Cal E}\text{\it{nd}}_{k}(\cE)~;~
Pf-fP\text{ is $\cO_{\P^n}$-linear }(f\in\cO_{\P^n})
\}\\
&\simeq \cE\ot D_{\P^n}^{\leq1}\ot \cE^*,
\endalign$$
which admits the exact sequence
$$
0\lra \cE\ot\cE^* \lra \Sigma'_{\cE}\lra \cE\ot\cE^*\ot T_{\P^n}\lra0.
$$

We define $\Sigma_{\cE}$ to be the inverse image of
$T_{\P^n}\simeq \cO_{\P^n}\ot T_{\P^n}\hra \cE\ot\cE^*\ot T_{\P^n}$
where $\cO_{\P^n}\hra \cE\ot\cE^*$ is the diagonal embedding. By definition
we have the exact sequence
$$
0\lra \cE\ot\cE^* \lra \Sigma_{\cE}\lra T_{\P^n}\lra0.
$$

\Lem 2-3. \it We have the isomorphism
$\pi_*\Sigma_{\cL}\os{\sim}{\lra}\Sigma_{\cE}.$
\rm
\demo{Proof}
It is easy to see that the image of the natural map
$\pi_*\Sigma_{\cL}\ra\Sigma'_{\cE}$ is contained in the sheaf
$\Sigma_{\cE}$. Since the sheaf $\pi_*\Sigma_{\cL}$ is generated
by the sections (2-4), it is surjective.
Due to the exact sequences
$$
0\lra \cO_{\P^n}\lra \pi_*\Sigma_{\cL}\lra \pi_*T_{\P}\lra0,
$$
$$
0\lra \pi_*T_{\P/\P^n}\lra \pi_*T_{\P}\lra T_{\P^n}\lra0,
$$
$$
0\lra \cO_{\P^n}\lra \cE\ot\cE^*\lra \pi_*T_{\P/\P^n}\lra0,
$$
we can see that $\pi_*\Sigma_{\cL}$ is a locally free sheaf
of rank $n+(r+s)^2$, which is the same as the one of $\Sigma_{\cE}$.
Thus the surjective map $\pi_*\Sigma_{\cL}\ra\Sigma_{\cE}$ is also
injective.
\qed
\enddemo

\Def 3-3. \it Define
$$
\Sigma_{\cE}^0=\pi_*\Sigma_{\cL}(-\log\P_*)\subset \Sigma_{\cE},
\qaq
\Sigma_{\cE}^{00}=\pi_*\Sigma_{\cL}(-\log\P_*+\Q_*)\subset
\Sigma_{\cE}.
$$
\rm\vskip 4pt

By (2-7) we have the commutative diagram of the exact sequences
$$
\CD
0@>>> \cE\ot\cE^* @>>> \Sigma_{\cE}@>>> T_{\P^n}@>>>0\\
@.@AAA@AAA@AA{=}A\\
0@>>>\cE\ot\cE_0^*\op \cO_{\P^n}^{\op s} @>>> \Sigma_{\cE}^0@>>>T_{\P^n}@>>>0\\
@.@AAA@AAA@AA{=}A\\
0@>>> \cO^{\op r+s}_{\P^n}@>>> \Sigma_{\cE}^{00}@>>> T_{\P^n}@>>>0,
\endCD
\leqno(2-8)$$
We have the global section
$$\sigma'=\pi_*\sigma=(F_i,G_j)_{1\leq i\leq r, 1\leq j\leq s}
\in H^0(\cE)={\us{1\leq i\leq r}{\op}}H^0(\cO(d_i))\op
{\us{1\leq j\leq s}{\op}}H^0(\cO(e_j)).$$
It induces the surjective map
$$j(\sigma'):\Sigma_{\cE}^0\lra \cE,
\quad P\longmapsto P(\sigma')
\tag{2-9}$$
that by (2-8) induces the exact sequence
$$0\lra T_{\P^n}(-\log X_*+Y_*) \lra \Sigma_{\cE}^{00}
@>{j(\sigma')}>> \cE \lra 0.
\tag{2-10}$$
where
$X_*=\sum_{1\leq j\leq s}X_i$
and
$Y_*=\sum_{1\leq j\leq s}Y_j$
with
$X_i:F_i=0 \subset \P^n$ and $Y_j:G_j=0 \subset \P^n$.


\vskip 20pt

\input amstex
\documentstyle{amsppt}
\hsize=16cm
\vsize=23cm

\head \S3. Proof of Theorem (I) \endhead
\vskip 8pt

In this section we complete the proof of Thoerem(I).
First we prove the following.

\Th 3-1. \it Let the assumption be as in Theorem(I).
For integers $q,\ell$ with $0 \leq q\leq n-r$ and $\ell\geq 0$
there is a natural isomorphism
$$ B_q(\d+\e-n-1+\ell) \isom
H^{q}(X,\WXZ {n-r-q}(\ell))_{prim}.$$
\rm\vskip 6pt

We start with the following lemma.

\Lem 3-1. \it
Let $q\geq 0$ and $\ell$ be integers.
Assuming $\l>-\d-\e$, there is a natural isomorphism
$$
\phi_q(\ell):A_q(\d+\e-n-1+\l)/J(\Sigma_{\cE}^{00}) \isom
H^q(\P^n,\Omega_{\P^n}^{n-q}(\log X_*+Y_*)(\l)).
$$
where
$J(\Sigma_{\cE}^{00})\subset A_q(\d+\e-n-1+\l)$
is the image of the map (cf. (2-9))
$$j(\sigma')\ot 1:H^0(\Sigma_{\cE}^{00}\ot S^{q-1}(\cE)\ot\cO(\l))
\ra H^0(S^q(\cE)\ot\cO(\l)))=A_q(\d+\e-n-1+\l).$$
\rm
\demo{Proof}
(2-10) gives rise to the following exact sequence
$$\multline
0\ra \Omega^{n-q}_{\P^n}(\log X_*+Y_*)(\l) \ra
\os{q}{\wedge}\Sigma_{\cE}^{00}\ot
\cO(\d+\e-n-1+\l)\ra \cdots\\
\ra \Sigma_{\cE}^{00}\ot
S^{q-1}(\cE)\ot\cO(\d+\e-n-1+\l)
\ra S^q(\cE)\ot\cO(\d+\e-n-1+\l) \ra 0.
\endmultline$$
Thus we have the natural map
$$
A_q(\d+\e-n-1+\l)/J(\Sigma_{\cE}^{00})\lra
H^q(\Omega_{\P^n}^{n-q}(\log X_*+Y_*)(\l)).
$$
In order to show that this is an isomorphism, it suffices
to show the following:
$$
H^b(\os{b}{\wedge}\Sigma_{\cE}^{00}\ot
S^{q-b}(\cE)\ot\cO(\d+\e-n-1+\l))=0, \quad
\text{for } 1\leq b\leq q,
$$
$$
H^{b-1}(\os{b}{\wedge}\Sigma_{\cE}^{00}\ot
S^{q-b}(\cE)\ot\cO(\d+\e-n-1+\l))=0, \quad
\text{for } 2\leq b\leq q.
$$
By the bottom sequence of (2-8) we have a decreasing filtration
$F$ on $\os{b}{\wedge}\Sigma_{\cE}^{00}$ such that
$$Gr^{\nu}_F(\os{b}{\wedge}\Sigma_{\cE}^{00})
=\bigoplus_{\binom {r+s}{\nu}} \os{b-\nu}{\wedge}T_{\P^n}
=\bigoplus_{\binom {r+s}{\nu}} \Omega_{\P^n}^{n-b+\nu}(n+1)
\quad(0\leq \nu\leq b). $$
Thus the above vanishing follows from the Bott
vanishing theorem.
\qed
\enddemo
\vskip 4pt
Now Th.(3-1) follows from the following two lemmas.
Recall the assumption (2-2).

\Lem 3-2. \it Assume $\ell\geq 0$.
Let $J(\Sigma_{\cL}(-\log\P_*))\subset A_q(\d+\e-n-1+\ell)$ be defined as
Lem.(2-2). Its image via $\phi_q(\ell)$ coincides with
the image of
$$
\us{1\leq \alpha\leq r}{\bigoplus}
H^q(\P^n,\Omega_{\P^n}^{n-q}(\log X_*^{(\alpha)}+Y_*)(\ell))
\ra
H^q(\P^n,\Omega_{\P^n}^{n-q}(\log X_*+Y_*)(\ell)).
$$
where
$X_*^{(\alpha)}=\sum_{\overset{i\not=\alpha}\to{1\leq i\leq r}}X_i$.
\rm

\Lem 3-3. \it
Suppose $n\geq r+1$ and put $X=X_1\cap \cdots \cap X_r$
and $Z_j=X\cap Y_j$.
By the assumption (2-2) $X\subset \P^n$ is a nonsingular complete
intersection of codimension $r$ and $\Zst=\Sigma_{1\leq j\leq s}Z_j$ is
a normal crossing divisor in $X$.
Let $r \leq a \leq n$ and $\ell\geq 0$.
Then the sequence
$$ \bigoplus_{1\leq \alpha\leq r}
H^{n-a}(\Omega^a_{\P^n}(\log X_*^{(\alpha)}+Y_*)(\ell))
\ra H^{n-a}(\Omega^a_{\P^n}(\log X_*+Y_*)(\ell))
\ra H^{n-a}(\WXZ {a-r}(\ell))_{prim} \ra 0$$
is exact where the last map is the composite of the successive residue maps
along $X_i$ $(1\leq i\leq r)$.
\rm

\vskip 6pt\noindent
\it Proof of Lem.(3-2) \rm
First we claim that we may show Lem.(3-2) replacing
$H^q(\Omega_{\P^n}^{n-q}(\log X_*^{(\alpha)}+Y_*)(\ell))$
with
$H^q(\Omega_{\P^n}^{n-q}(\log X_*+Y_*)(-X_\alpha)(\ell))$.
The claim follows from the general lemma (3-4) below.
The exact sequence
$$\multline
0\ra \Omega^{n-q}_{\P^n}(\log X_*+Y_*)(-X_\alpha)(\ell) \ra
\os{q}{\wedge}\Sigma_{\cE}^{00}\ot\cO(\d+\e-n-1)(-X_\alpha)(\ell)\ra \cdots\\
\ra \Sigma_{\cE}^{00}\ot S^{q-1}(\cE)\ot\cO(\d+\e-n-1)(-X_\alpha)(\ell)
\ra S^q(\cE)\ot\cO(\d+\e-n-1)(-X_\alpha)(\ell) \ra 0
\endmultline$$
gives rise to a commutative diagram
$$
\CD
H^0(S^q(\cE)\ot\cO(\d+\e-n-1+\ell)(-X_\alpha))@>{\tau'}>>
H^q(\Omega_{\P^n}^{n-q}(\log X_*+Y_*)(-X_\alpha)(\ell))\\
@V{F_\alpha}VV@VVV\\
A_q(\d+\e-n-1+\ell)@>{\tau}>>
H^q(\Omega_{\P^n}^{n-q}(\log X_*+Y_*)(\ell)).
\endCD
$$
Thus it suffices to show
$$J(\Sigma_{\cL}(-\log\P_*))=J(\Sigma_{\cE}^{00})+\sum_{1\leq \alpha\leq r}
\Im(F_\alpha)\subset A_q(\d+\e-n-1+\ell)\leqno(*)$$
and that $\tau'$ is surjective.
By (2-7) $(*)$ holds if we replace $J(\Sigma_{\cE}^{00})$ by
$$
J(\Sigma_{\cL}(-\log P_*+Q_*)):=\Im\big(
H^0(\Sigma_{\cL}(-\log\P_*+\Q)\ot \cL^{q-1}\ot\pi^*\cO(\l'))
@>{j(\sigma)\ot 1}>> H^0(\cL^q\ot\pi^*\cO(\l'))\big)
$$
with $\ell'=\d+\e-n-1+\ell$. Therefore it suffices to show
$J(\Sigma_{\cE}^{00})=J(\Sigma_{\cL}(-\log \P_*+\Q_*))$
that follows from the isomorphism
$$
\pi_*\Sigma_{\cL}(-\log\P_*+\Q_*)\ot
\pi_*\cL^q\os{\sim}{\ra}\pi_*(\Sigma_{\cL}(-\log\P_*+\Q_*)\ot\cL^q),
$$
which follows from (2-7).
To show the surjectivity of $\tau'$, it suffices to show that
$$H^b(\os{b}{\wedge}\Sigma_{\cE}^{00}\ot
S^{q-b}(\cE)\ot\cO(\d+\e-n-1)(-X_\alpha)(\ell))=0 \quad
\text{for } 1\leq b\leq q.
$$
This follows from the Bott vanishing theorem by the same argument as
the proof of Lem.(3-1) using (2-8).
\qed

\vskip 6pt\noindent
\it Proof of Lem.(3-3) \rm
We can directly check the exactness of the following sequence
$$\multline
0\ra \Omega^a_{\P^n}(\log Y_*)\ra\cdots
\ra \us{1\leq \alpha<\beta\leq r}{\op}
\Omega^a_{\P^n}(\log X_*^{(\alpha\beta)}+Y_*) \ra \\
\ra \us{1\leq \alpha\leq r}{\op}\Omega^a_{\P^n}(\log X_*^{(\alpha)}+Y_*)
\ra \Omega^a_{\P^n}(\log X_*+Y_*)
\ra \WXZ {a-r}\ra 0,
\endmultline
\tag{3-1}
$$
where
$X_*^{(\alpha\beta)}=\sum_{\overset{i\not=\alpha,\beta}\to{1\leq i\leq s}}X_i$
and so on.
Thus the desired assertion follows from the following general lemma.

\Lem 3-4. \it
Let the notaion and the assumption be as Def.(1-2) and let $d=\dim(X)$.
\roster
\item
If $a+b\not=\dim X$ and $a\geq 1$ and $\ell\geq 0$,
$ H^a(X,\WXZ b(\ell))_{prim}=0$.
\item
For $1\leq \alpha\leq r$ and $a\geq 1$ the natural map
$$H^a(X,\WXZ {d-a}(-Z_\alpha)(\ell))\to H^a(X,\WXZa {d-a}(\ell))_{prim}
$$
is surjective where
$\Zast=\underset{\overset{j\not=\alpha}\to{1\leq j\leq s}}\to{\sum} Z_j$.
\endroster
\rm
\demo{Proof}
First we note that $ H^a(X,\Omega_{X}^{b}(\ell))_{prim}=0$ if $a+b\not=\dim X$
and$a\geq 1$ and $\ell\geq 0$. This is well-known fact on cohomoology of
smooth complete intersections. We also note that the case $b=0$ of Lem.(3-4)(1)
is true by the same reason. Assuming $b\geq 1$ we have the exact sequence
$$
H^a(X,\Omega_{X}^{b}(\log \Zst^{(1)})(\ell)) \to
H^a(X,\WXZ {b}(\ell)) \to
H^a(Z_1,\Omega_{Z_1}^{b-1}(\log \Wst^{(1)})(\ell))
\quad (\Wst^{(1)}=\sum_{2\leq j\leq s} Z_1\cap Z_j)
$$
By induction on $s$ and $b$ we are reduced to the case $\dim(X)=1$.
In this case we have only to consider
$  H^1(X,\WXZ {1}(\ell))_{prim}$
which we easily see vanishes.
This completes the proof of Lem.(3-4)(1).
Lem.(3-4)(2) follows from (1) in view of the exact sequence
$$ 0\to \WXZ b(-Z_\alpha)\to \WXZa {b} \to
\Omega_{Z_\alpha}^{b}(\log \Wst^{(\alpha)}) \to 0
\quad
(\Wst^{(\alpha)}=
\underset{\overset{j\not=\alpha}\to{1\leq j\leq s}}\to{\sum} Z_\alpha\cap Z_j)
$$
\qed
\enddemo
\vskip 8pt

Next we show Theorem(I)(2), namely the follwoing statement.
\Lem 3-5. \it Let the assumption be as Lem.(3-2).
There is a natural map
$$ B_1(0) \lra H^1(X,\TXZ)_{alg}$$
which is an isomorphism if $\dim(X)\geq 2$.
\rm
\demo{Proof}
The section $(F_i,G_j)_{1\leq i\leq r,1\leq j\leq s} \in H^0(\cE) $
(cf. (2-1)) defines the surjective map
$$ j_1:\cE_0^* \ot \cE \lra I_X \ot \cE, \quad \xi^*_i\ot \cdot \longmapsto f_i \ot \cdot
$$
($I_X$ denotes the ideal sheaf of $X$) and the map
(We denote $\Cal O=\cO_{\P^n}$)
$$
 j_2: \Cal O^{\oplus s} \lra \cE,\quad
 e_j=(0,\cdots,1,\cdots,0) \longmapsto g_j\eta_j\quad(1 \leq j \leq s)$$
Here $\xi_i$ (resp. $\eta_j$) is a local frame of $\Cal O(d_i)$
(resp. $\Cal O(e_j)$) and
$\Sigma f_i \xi_i+\Sigma g_j\eta_j$ is the local description of the image
of $\sum^r_{i=1} F_i\mu_i+\sum^s_{j=1}G_j\ld_j \in H^0(\PE,\cL)$
under the isomorphism
$$H^0(\P, \cL) \simeq
H^0(\P^n, \Cal E)=
\os{r}{\us{i=1}{\op}}H^0(\P^n, \Cal O(d_i))\oplus
\os{s}{\us{j=1}{\op}}H^0(\P^n, \Cal O(e_j)).$$
Put
$$
I=\Im(j_1+j_2: (\cE_0^* \ot \cE)\op \Cal O^{\oplus s}  \lra \cE),
$$
which is generated by local sections
$$
f_i \xi_{i'}, \quad  f_i\eta_j, \quad  g_j\eta_j \qquad (1 \leq i,i' \leq r, \quad 1 \leq j \leq s).
$$
We have the following commutative diagram (cf. (2-8) and (2-9))
$$
\matrix
       & 0 &
       & 0 &
       & 0 &
\\
       & \downarrow &
       & \downarrow &
       & \downarrow &
\\
0 \lra & L &
  \lra & K &
  \lra & T &
  \lra 0
\\
       & \downarrow &
       & \downarrow &
       & \downarrow &
\\
0 \lra & \Cal E_0^* \ot \Cal E \op \Cal O^{\oplus s}&
  \lra & \Sigma_{\cE}^{0} &
  \lra & T_{\Bbb P^n} &
  \lra 0
\\
       & \qquad \downarrow j_1 + j_2 &
       & \qquad \downarrow j(\sigma')&
       & \quad \downarrow j_3 &
\\
0 \lra & I&
  \lra & \Cal E &
  \lra & \Cal E/I&
  \lra 0
\\
       & \downarrow &
       & \downarrow &
       & \downarrow &
\\
       & 0 &
       & 0 &
       & 0 &
\endmatrix
\tag 3-2
$$
The map $j_3$ in the above diagram can be written as follows:
$$
j_3: \frac{\partial}{\partial x} \longmapsto
\sum^r_{i=1}\frac{\partial f_i}{\partial x}\xi_i +
\sum^s_{j=1}\frac{\partial g_j}{\partial x}\eta_j  \mod I
$$
and it is easy to see that this implies the following exact sequence
$$
0 \lra I_X \ot T_{\P^n} \lra T \lra \TXZ \lra 0.
\tag 3-3
$$
We get the map
$$
\align
B_1(0)
&=\Coker
  (H^0(\Bbb P,\Sigma_{\cL}(-\log \P_*) \lra H^0(\Bbb P,\Cal L))
\\
&\simeq \Coker
  (H^0(\Bbb P^n,\Sigma_{\cE}^{0}) \lra H^0(\Bbb P^n,\Cal E))
\\
&\os{a}{\lra}
  H^1(\Bbb P^n,K)
   \text{ (from the middle vertical sequence in (3-2))}
\\
&\os{b}{\lra}
  H^1(\Bbb P^n,T)
   \text{ (from the top horizontal sequence in (3-2))}
\\
&\os{c}{\lra}
  H^1(X,\TXZ)
   \text{ (from (3-3))}
\endalign
$$
Thus Lem.(3-5) follows from the following.

\Lem 3-6. \it Assume $\dim(X)=n-r\geq 2$ and $n\geq 3$.
\roster
\item
$H^1(\Sigma_{\cE}^{0})=0.$
\item
$H^1(L)=H^2(L)=0.$
\item
$H^1(I_X \ot T_{\Bbb P^n})=0.$
\item
$\Ker(H^1(\TXZ) @>\delta>> H^1(I_X \ot T_{\Bbb P^n}))= H^1(\TXZ)_{alg}$
where $\delta$ is induced by (3-3).
\endroster
\rm
\demo{Proof}
(1) follows from (2-8) and the Bott vanishing.
To show (2) we consider the following commutative diagram
$$
\matrix
       & 0 &
       & 0 &
       & 0 &
\\
       & \downarrow &
       & \downarrow &
       & \downarrow &
\\
0 \lra & L_1 &
  \lra & L &
  \lra & L_2 &
  \lra 0
\\
       & \downarrow &
       & \downarrow &
       & \downarrow &
\\
0 \lra & \cE^*_0\ot\cE &
  \lra & \cE^*_0\ot\cE \op \Cal O^{\oplus s}  &
  \lra & \Cal O^{\oplus s} &
  \lra 0
\\
       & \qquad \downarrow j_1  &
       & \qquad \downarrow j_1+j_2&
       & \quad \downarrow j_2' &
\\
0 \lra & I_X \ot \cE&
  \lra & \Cal E &
  \lra & \Cal E \ot \Cal O_X&
  \lra 0
\\
       & \downarrow &
       &  &
       & &
\\
       & 0 &
       & &
       &  &
\endmatrix
$$
where $j_2': e_k \longmapsto g_k\eta_k \mod I_X$ ($1 \leq k \leq s$).
Therefore we have
$L_2=\Ker(j_2')=I_X^{\op s}.$
From the Koszul exact sequence
$$
0 \lra \os{r}{\wedge} \Cal E_0^*
      \lra \cdots
        \lra \os{2}{\wedge} \Cal E_0^*
         \lra \Cal E_0^*
          \lra I_X
            \lra 0,
$$
we can see that $L_i$ has the following resolution.
$$
\aligned
&0 \lra (\os{r}{\wedge} \Cal E_0^*)^{\op s}
      \lra \cdots
        \lra \Cal E_0^{* \op s}
\lra L_2
            \lra 0,
\\
&0 \lra \os{r}{\wedge} \Cal E_0^* \ot \cE
      \lra \cdots
        \lra \os{2}{\wedge} \Cal E_0^*\ot \cE
          \lra L_1
            \lra 0.
\endaligned
$$
Therefore (2) follows from the Bott vanishing.
(3) is an easy consequence of the Euler exact sequence
$$ 0\to \cO_{\P^n} \to  \cO_{\P^n}(1)^{\oplus n+1} \to T_{\P^n} \to 0.
\leqno(3-4)$$
Finally we see easily that (4) is reduced to prove
$$\Ker(H^1(X,\TX) @>\delta_1>> H^2(I_X \ot T_{\Bbb P^n}))=
\Ker(H^1(X,\TX) @>\delta_2>> H^2(X,\cO_X)),\leqno(*)$$
where $\delta_1$ is induced by the exact sequence
$$ 0 \lra I_X \ot T_{\P^n} \lra T' \lra \TX \lra 0
\qwith
T'=\Ker(T_{\P^n}\to \cE_0\otimes \cO_X)$$
and $\delta_2$ is the map in Def.(1-3). By the commutative diagram
$$  \matrix
&&& 0&& 0\\
&&&\downarrow&&\downarrow\\
0 \lra  & I_X \ot T_{\P^n} & \lra & T' & \lra & \TX & \lra 0\\
& \| &&\downarrow&&\downarrow\rlap{$\iota$}\\
0 \lra & I_X \ot T_{\P^n} & \lra & T_{\P^n} & \lra & T_{\P^n}\ot\cO_X &\lra 0\\
&&&\downarrow&&\downarrow\\
&&& \cE_0\otimes \cO_X &=&\cE_0\otimes \cO_X\\
&&&\downarrow&&\downarrow\\
&&& 0&& 0\\
\endmatrix$$
and the fact $H^2(T_{\P^n})=0$, we may prove $(*)$ after replacing
the left hand side with the kernel of
$H^1(X,\TX) @>\iota>> H^1(T_{\Bbb P^n}\ot\cO_X).$
(3-4) induces the boundary map
$H^1(T_{\Bbb P^n}\ot\cO_X) @>\delta_3>> H^2(\cO_X)$ which is injective.
We see that
the composite of $\iota$ and $\delta_3$ coincides with the map in Def.(1-3)
by noting that the extension class of (3-4) is given by
$c_1(\cO_{\P^n}(1))\in H^1(\Omega^1_{\P^n})$. This completes the proof.
\qed
\enddemo
\enddemo


\vskip 20pt

\input amstex
\documentstyle{amsppt}
\hsize=16cm
\vsize=23cm

\head \S4. A vanishing lemma \endhead
\vskip 8pt

The following result is the technical heart of the proof of
Theorem(II) and (III). For a vector bundle $\cF$, $\cF^*$ denote its dual.
\vskip 5pt

\Th 4-1)(vanishing lemma. \it Assume $s\geq 1$.
Let $p$, $w$, $\nu$, $\l$ be integers. Then
$$
H^w(\P, \os{p}{\wedge}\Sigma_{\cL}(-\log \P_*)^* \ot \L^{\nu}\ot\pi^*\Cal O(\l))=0
$$
if one of the following conditions is satisfied.
We put $m=n+r+s-1=\dim\P$.
\roster
\item
$w>0$, $\nu \geq -s+1$, $\l \geq 0$ and $(\nu,\l) \not=(0,0)$
\item
$p-\nu \leq w < m$ and $\nu \leq -1$
\item
$e_1=e_2=\cdots=e_s$ and $0<w\not=n$ and $\nu\geq -s+1$
\item
$e_1=e_2=\cdots=e_s$, $0<w$, $\nu=\ell=0$ and $p\leq n$
\item
$p-\nu\leq r+s-1$ and $\nu\leq -1$.
\endroster
\roster
\item"$(1)^*$"
$w<m$, $\nu \leq -1$, $\l \leq \e$ and $(\nu,\l) \not=(-s,\e)$
\item"$(2)^*$"
$0 < w < p-\nu-s$ and $\nu \geq -s+1$
\item"$(3)^*$"
$e_1=e_2=\cdots=e_s$ and  $m>w\not=r+s-1$ and $\nu<0$
\item"$(4)^*$"
$e_1=e_2=\cdots=e_s$, $m>w$, $\nu=-s$, $\ell=\e$ and $p\geq r+s$
\item"$(5)^*$"
$p-\nu\geq n+1$ and $\nu\geq -s+1$.
\endroster
\rm

\demo{Proof}
For simplicity we put $\Sigma=\Sigma_{\cL}(-\log \P_*)$.
By the Euler exact sequence (cf. Lem.(2-1)(3))
we have the isomorphism
$$
\Omega^1_{\P/\P^n}(\log \P_*) \simeq
(\pi^*\cE_0 \ot \cL^{-1})\op \Cal O_{\P}^{\op s-1}.
\tag{4-1}$$
and we have
$$
\os{m+1}{\wedge}\SS = \text{ det }\Omega^1_{\P}(\log \P_*)
=\cL^{-r}\ot\pi^*\Cal O(\d-n-1).
\tag{4-2}
$$
Noting
$$ K_{\P}=\text{ det }\Omega^1_{\P}=
\pi^*(K_{\P^n}\otimes \text{det}{\Cal E})\ot\cL^{-r-s}
=\cL^{-r-s}\ot\pi^*{\Cal O}(\d+\e-n-1),
\tag{4-3}
$$
we have the Serre duality:
$$
H^w(\os{p}{\wedge}\SS \ot \cL^{\nu}\ot\pi^*\Cal O(\l))^*\simeq
H^{m-w}(\os{m+1-p}{\wedge}\SS \ot \cL^{-\nu-s}\ot\pi^*\Cal O(\e-\l)).
$$
Therefore the assertion ($n$) ($1 \leq n \leq 5$) is equivalent to $(n)^*$ and
we only need to show (1),(2),(3),(4) and (5).
By Def.(2-1) we have the exact sequence
$$0\to\cO_{\P}\to \Sigma \to T_{\P}(-\log\P_*)\to 0$$
that induces the exact sequence
$$
0 \ra \Omega^p_{\P}(\log \P_*)
\ra \os{p}{\wedge}\SS \ra
\Omega^{p-1}_{\P}(\log \P_*)
\ra 0.
\tag{4-4}
$$
Moreover the exact sequence
$$0 \ra \pi^*\Omega^1_{\P^n} \ra \Omega^1_{\P}(\log \P_*) \ra
\Omega^1_{\P/\P^n}(\log \P_*) \ra 0$$
gives rise to a finite decreasing filtration $F^{\cdot}$
on $\Omega^q_{\P}(\log \P_*)$ such that
$$
\text{Gr}^{a}_F(\Omega^q_{\P}(\log \P_*))
=\pi^*\Omega^{a}_{\P^n} \ot \Omega^{q-a}_{\P/\P^n}(\log \P_*)\simeq
\os{q-a}{\us{i=0}{\op}} \pi^*\Omega^a_{\P^n}\ot
[\os{i}{\wedge}\pi^*\cE_0\ot \cL^{-i}]^{\oplus\binom{s-1}{q-a-i}}
\tag{4-5}
$$
where the second isomorphism follows from (4-1).
Hence we obtain the spectral sequence
\def\sE#1#2#3#4{{_{#4}E_{#1}^{#2,#3}}}
$$\sE 1 a b q=
\os{q-a}{\us{i=0}{\op}}
H^{a+b}(\P,\cL^{\nu-i}\ot\pi^*(\Omega^{a}_{\P^n}(\l)\ot\os{i}{\wedge}\cE_0))^{\oplus\binom{s-1}{q-a-i}}
\Rightarrow
H^{a+b}(\P,\Omega^{q}_{\P}(\log \P_*)\ot\cL^{\nu}\ot\pi^*\Cal O(\ell))
\tag{4-6}
$$
Noting that $\os{i}{\wedge}\cE_0=0$ for $i>r$, Lem.(2-1)(2) implies
$$
\eqalign{
H^w(\P, \cL^{\nu-i}\ot&\pi^*(\Omega^{a}_{\P^n}(\l)\ot\os{i}{\wedge}\cE_0))\cr
& \simeq \left \{\eqalign{
&H^w(\P^n, S^{\nu-i}(\cE)\ot \Omega^{a}_{\P^n}(\l)\ot\os{i}{\wedge}\cE_0)
\quad \hbox{ if } \nu\geq i, \cr
&H^{w-t+1}(\P^n, S^{-\nu+i-t}(\cE^*)\ot \Omega^{a}_{\P^n}\ot\os{i}{\wedge}\cE_0(\l-\d-\e))
\quad \hbox{ if } \nu\leq i-t, \cr
&\hskip 10pt 0 \qquad \hbox{ if } i>r \hbox{ or }i-t<\nu<i. \cr
}\right.\cr}
\leqno(4-7)$$
Here we put $t=r+s$. By (4-6) and (4-4) the desired vanishing
in cases (1), (2) and (5) follows from the following.

\proclaim{{\it Claim}}
$(i)$  Under the assumption (1), $\sE 1 a b q=0$ if $w=a+b>0$.
\vskip 4pt\noindent
$(ii)$  Under the assumption (2), $\sE 1 a b q=0$ if $w=a+b>0$ and $q\leq p$.
\vskip 4pt\noindent
$(iii)$  Under the assumption (5), $\sE 1 a b q=0$ if $q\leq p$.
\endproclaim

First assume $\nu\geq -s+1$. For $i\leq r$, $\nu\geq i-t+1$.
Hence the first assertion of the claim follows immediately from (4-7)
and the Bott vanishing. Next assume (2) and $q\leq p$.
We have $\nu-i \leq \nu \leq -1$ so that by (4-7) we may assume
$\nu-i \leq -t$. Hence $a\leq q-i \leq p-i \leq p-\nu-t$.
By (4-7) and the Bott vanishing
we get $\sE 1 a b q=0$ if $p-\nu-t < w-t+1 <n$, that is,
$p-\nu \leq w < m$. This completes the proof in case (2).
Finally the assertion in case (5) follows from (4-7) by noting
$0\leq i\leq q-a$ in (4-6).
\par

Next we treat the (3). Assume $\nu\geq -s+1$. By (4-7)
$$ \sE 1 a b q=
\os{q-a}{\us{i=0}{\op}}
H^{a+b}(\P^n, S^{\nu-i}(\cE)\ot \Omega^{a}_{\P^n}(\l)\ot\os{i}{\wedge}\cE_0)
^{\oplus\binom{s-1}{q-a-i}}
\tag{4-8}
$$
By the Bott vanishing,
$\sE 1 a 0 q,\sp\sE 1 a {-a} q,\sp \sE 1 a {n-a} q$ with
$0\leq a\leq \min\{n,q\}$
are the only terms which are possibly non-zero. Therefore we have

\roster
\item"$(i)$"
$\sE 2 a b q =\sE {\infty} a b q$ for any $a,b$.
\item"$(ii)$"
$\sE 1 a {-a} q=\sE {\infty} a {-a} q$ and
$\sE 1 a {n-a} q =\sE {\infty} a {n-a} q$ for $a\not=0,n$.
\item"$(iii)$"
$\sE 2 a 0 q \simeq
H^{a}(\P,\Omega^{q}_{\P}(\log \P_*)\ot\cL^{\nu}\ot\pi^*\Cal O(\ell))
\quad\hbox{ for } a\not=0,n.$
\endroster

Note that the boundary map coming from (4-4)
$$
\partial_a \sp:\sp
H^{a-1}(\P,\Omega^{p-1}_{\P}(\log \P_*)\ot\cL^{\nu}\ot\pi^*\Cal O(\ell))
\to
H^{a}(\P,\Omega^{p}_{\P}(\log \P_*)\ot\cL^{\nu}\ot\pi^*\Cal O(\ell))
$$
is induced by the cup product with the class
$$ \tilde c:= c_1(\cL)_{|\P\backslash \P_*} \in
H^1(\P,\Omega^{1}_{\P}(\log \P_*))$$

\proclaim{{\it Claim}}
Assume $e_1=e_2=\cdots=e_s=e$. Then the natural map
$$ \pi^*\sp:\sp H^1(\P^n,\Omega^1_{\P^n}) \to
H^1(\P,\Omega^{1}_{\P}(\log \P_*))$$
injective and we have
$\tilde c= \pi^*(c_1(\Cal O_{\P^n}(e))).$
\endproclaim

We know that
$H^1(\P,\Omega^{1}_{\P})$ is the direct sum of
$\pi^*(H^1(\P^n,\Omega^1_{\P^n}))$ and the subspace spanned by $c_1(\cL)$.
The kernel of
$H^1(\P,\Omega^{1}_{\P})\to H^1(\P,\Omega^{1}_{\P}(\log \P_*))$
is generated by $c_1(\Cal O_{\P}(\P_j))$ for $1\leq j\leq s$.
The claim follows from the fact
$\cL\ot\pi^*\Cal O(-e_j)= \Cal O_{\P}(\P_j)$.
\par

From the claim the map $\partial_a$ and the maps
$\delta \sp:\sp \sE 1 {a-1} b {p-1} @>{\cup c}>> \sE 1 a b p$
with $c$ a non-zero multiple of $c_1(\Cal O_{\P^n}(e))$,
are compatible with respect to the spectral sequence (4-6).
Thus we get the commutative diagram
$$
\matrix
&\sE 1 0 0 {p-1}&\to
&\sE 1 1 0 {p-1}&\to
&\cdots&\to
&\sE 1 {n-2}  0 {p-1}&\to
&\sE 1 {n-1}  0 {p-1}&\to
&\sE 1 n 0 {p-1}
\\
&&\searrow
&&\searrow
&&
&&\searrow
&&\searrow
\\
&\sE 1 0 0 {p}&\to
&\sE 1 1 0 {p}&\to
&\sE 1 2 0 {p}&\to
&\cdots&\to
&\sE 1 {n-1} 0 {p}&\to
&\sE 1 n 0 {p}
\\
\endmatrix
$$
where the horizontal arrows are the differentials in (4-6) and
the slanting arrows are the above map $\delta$. To show the desired vanishing
in case (3), we need show that if $w\not=0,n$ then $\partial_{w}$ is
surjective and
$\partial_{w+1}$ is injective. Therfore it suffices to show the following

\proclaim{{\it Claim}}
The map
$\delta \sp:\sp \sE 1 {a-1} 0 {p-1} @>{\cup c}>> \sE 1 a 0 p$
is an isomrophism for $2\leq a\leq n-1$ and surjective for $a=1$ and
injecctive for $a=n$.
\endproclaim

By (4-8), for $a\not=0,n$, we have the isomorphism
$$
\sE 1 a 0 q \simeq
H^a(\P^n, \Omega^{a}_{\P^n})^{\oplus \phi}.
$$
Here $\phi$ is determined as follows. If $q<a$, $\phi=0$.
If $q\geq a$, writing
$$
\os{q-a}{\us{i=0}{\op}}
[S^{\nu-i}(\cE)\ot \os{i}{\wedge}\cE_0\ot
\Cal O(\ell)]^{\oplus\binom{s-1}{q-a-i}}
= {\us{k}{\op}} \Cal O(\ell_k),$$
we put $\phi=\# \{k|\ell_k=0\}$. Thus, fixing $\nu$ and $\ell$,
$\phi=\phi(q-a)$ is a function of $q-a$. Note that even for $a=0$ or $n$,
we have the injection
$ H^a(\P^n, \Omega^{a}_{\P^n})^{\oplus \phi} \hookrightarrow \sE 1 a 0 q$.
Since
$ H^{a-1}(\P^n, \Omega^{a-1}_{\P^n}) @>{\cup c}>>H^{a}(\P^n, \Omega^{a}_{\P^n})$
is clearly an isomorphism, the claim is proven.

Finally we treat the case (4). By the vanishing lemma(3) it suffices to show
$H^n(\P,\os{p}{\wedge}\SS)=0$. By (4-4) this is reduced to proving the
following.
\roster
\item"$(i)$"
$H^n(\P,\Omega^{p-1}_{\P}(\log \P_*))=0.$
\item"$(ii)$"
$H^{n-1}(\P,\Omega^{p-1}_{\P}(\log \P_*))=0
@>{\cup \tilde c}>> H^n(\P,\Omega^{p}_{\P}(\log \P_*))$ is surjective.
\endroster

As before we have the spectral sequence
$$\sE 1 a b q
\Rightarrow
H^{a+b}(\P,\Omega^{q}_{\P}(\log \P_*))
$$
where
$$\sE 1 a b q
 \simeq \left \{\eqalign{
&H^{a+b}(\P^n, \Omega^{a}_{\P^n})^{\oplus\binom{s-1}{q-a}}
\quad \hbox{ if } 0\leq a\leq \min\{q,n\}, \cr
&\hskip 10pt 0 \qquad \hbox{ otherwise }. \cr
}\right.
$$
Thus $\sE 1 a b q =0$ unless $b=0$ and the spectral sequence degenerates at
$E_2$. Now the assertion $(i)$ follows from the fact that $\sE 1 n 0 {p-1} =0$
by the assumption $p\leq n$. The assertion $(ii)$ can be shown by the same
argument as the proof in case (3).
\qed
\enddemo


\vskip 20pt

\input amstex
\documentstyle{amsppt}
\hsize=16cm
\vsize=23cm

\head \S5. Proof of Theorem(II) \endhead
\vskip 8pt

In this section we prove Thoerem(II).
We deduce it from the Serre duality theorem and Lem.(4-1).
The exact sequence (cf.(2-6))
$$ 0 @>>> T_{\P}(-\log \cZ+\P_*) @>>> \Sigma_{\cL}(-\log\P_*)
@>{j(\sigma)}>> \cL @>>> 0,$$
and (4-2) induce the Koszul exact sequence
$$
0 \ra \L^{-m-1} \ra \SS\ot\L^{-m}\ra
\cdots \ra \os{m+1}{\wedge}\SS \ra 0.\quad (\Sigma=\Sigma_{\cL}(-\log\P_*))
\tag 5-1
$$
Tensoring with $\L^{r+q}\ot \pi^*\Cal O(\l)$, we get the exact sequence
where we denote $\cO=\cO_{\bP^n}$
$$
0 \ra \L^{-m+r+q-1}\ot\pi^*\Cal O(\l) \ra
\cdots \ra \os{m}{\wedge}\SS \ot \L^{r+q-1}\ot\pi^*\Cal O(\l)
@>\delta>>
\os{m+1}{\wedge}\SS \ot \L^{r+q}\ot\pi^*\Cal O(\l)\ra 0.
\tag{5-2}
$$
By (4-2) we have
$$\os{m}{\wedge}\SS \ot \L^{r+q-1}\ot\pi^*\Cal O(\l)
=\Sigma \ot \L^{q-1}\ot\pi^*\Cal O(\d-n-1+\l),$$
$$\os{m+1}{\wedge}\SS \ot \L^{r+q}\ot\pi^*\Cal O(\l)
=\L^{q}\ot\pi^*\Cal O(\d-n-1+\l)$$
and the map $\delta$ in (5-2) is nothing but $j(\sigma)\ot 1$.
Therefore we have the canonical map
$$
B_q(\d-n-1+\l)
\ra \Ker[ H^m(\L^{-m+r+q-1}\ot \pi^*\Cal O(\l))\ra H^m(\SS\ot \L^{-m+r+q}\ot \pi^*\Cal O(\l))]
$$
By (4-3) and the Serre duality the right hand side is isomorphic to the dual
of
$$
\multline
\Coker[H^0(\Sigma \ot \L^{n-r-q-1}\ot\pi^*\Cal O(\d+\e-n-1-\l))
\ra H^0(\L^{n-r-q}\ot\pi^*\Cal O(\d+\e-n-1-\l))]\\
=B_{n-r-q}(\d+\e-n-1-\l).
\endmultline
$$
Thus we get the canonical maps
\def\h#1#2{h_{#1}(#2)}
\def\hst#1#2{h^*_{#1}(#2)}
$$ \aligned
&\h q \ell\sp:\sp B_q(\d-n-1+\l) \ra B_{n-r-q}(\d+\e-n-1-\l)^*,\\
&\hst q \ell\sp:\sp B_{n-r-q}(\d+\e-n-1-\l) \to B_q(\d-n-1+\l)^*,\\
\endaligned$$
where $\hst q \ell$ is the dual of $\h q \ell$.
In particular we get the \it trace map \rm
$$ \tau:=\h {n-r} {\d+\e-n-1}\sp:\sp
B_{n-r}(2(\d-n-1)+\e) \ra B_{0}(0)^*=k.
\tag{5-3}
$$
For $x\in B_q(\d-n-1+\l)$ and $y\in B_{n-r-q}(\d+\e-n-1-\l)$ we let
$<x,y>\in k$ be the evaluation of $\h q \ell(x)$ at $y$.
This gives us a bilinear pairing
$$<\sp,\sp>\sp:\sp B_q(\d-n-1+\l) \otimes B_{n-r-q}(\d+\e-n-1-\l) \to k
\sp;\sp (x,y)\to <x,y>.$$

\Lem 5-1. \it
(1) $\hst q \ell$ coincides with
$$\h {n-r-q} {\e-\ell} \sp:\sp B_{n-r-q}(\d+\e-n-1-\ell) \to
B_{q}(\d-n-1+\ell)^*.$$
\vskip 4pt\noindent
(2) We have
$$<x,y>=\tau(xy)\qfor x\in B_q(\d-n-1+\l) \hbox{ and }
y\in B_{n-r-q}(\d+\e-n-1-\l),$$
where $xy\in B_{n-r}(2(\d-n-1)+\e)$ is the multiplication of $x$and $y$.
\rm \demo{Proof}
The first assertion follows by taking the dual of (5-2) in view of the
isomorphism
$$\big(\os{p}\wedge \SS\otimes\L^\nu\otimes\pi^*\Cal O(\ell)\big)^*\ot K_\P
\simeq \os{m+1-p}\wedge \SS\otimes\L^{-\nu-s}\otimes\pi^*\Cal O(\e-\ell).$$
The second assertion follows from the commutative diagram
$$
\matrix
0 \ra &\L^{-m+r+q-1}\ot\pi^*\Cal O(\l)& \ra  \cdots \hskip 200pt\hbox{ }\\
&\downarrow\\
0 \ra &\L^{-r-s}\ot\pi^*\Cal O(\d+\e-n-1)& \ra  \cdots\hskip 200pt\hbox{ }\\
\endmatrix
$$
$$
\matrix
\cdots\ra &\os{m}{\wedge}\SS \ot \L^{r+q-1}\ot\pi^*\Cal O(\l)&\ra&
\os{m+1}{\wedge}\SS \ot \L^{r+q}\ot\pi^*\Cal O(\l)&\ra 0\\
&\downarrow&&\downarrow\\
\cdots\ra &\os{m}{\wedge}\SS \ot \L^{n-1}\ot\pi^*\Cal O(\d+\e-n-1)&\ra&
\os{m+1}{\wedge}\SS \ot \L^{n}\ot\pi^*\Cal O(\d+\e-n-1)&\ra 0\\
\endmatrix
$$
where the cup product with an element of
$H^0(\L^{n-r-q}\otimes \pi^*\Cal O(\d+\e-n-1-\ell))$ gives the vertical maps.
\enddemo

\Lem 5-2. \it Assume $s\geq 1$ and $e_1=\cdots=e_s$.
\vskip 4pt\noindent
(1)  $\h q \ell$ is injective under one of the following conditions.
\roster
\item"$(i)$"  $n-r\geq q$ and $\ell<\e$.
\item"$(ii)$"  $n-r\geq q\geq \frac{n-r+1}2$.
\endroster
\vskip 4pt\noindent
(2)  $\h q \ell$ is an isomorphism under one of the following conditions.
\roster
\item"$(i)$"  $q<n-r$ and $\ell<\e$.
\item"$(ii)$"  $0<q<n-r$.
\item"$(iii)$"  $0\leq \ell\leq \e$ and $r+s\leq n$.
\endroster
\vskip 4pt\noindent
(3) Assuming $n-r\geq 1$, $\tau$ is an isomorphism.

\demo{Proof}
By the exact sequence (5-2), $\h q \ell$ is surjective if
$$
\text{(a) } H^a(\os{m+1-a}{\wedge}\SS \ot \L^{r+q-a} \ot \pi^*\Cal O(\l))=0
\text{ for } 1 \leq a \leq m-1,
$$
and is injective if
$$
\text{(b) } H^b(\os{m-b}{\wedge}\SS \ot \L^{r+q-b-1} \ot \pi^*\Cal O(\l))=0
\text{ for }1 \leq b \leq m-1.
$$

To show the injectivity, we show (b). First we assume (1)$(i)$.
By the assumption $\ell<\e$ and Lem.(4-1)$(1)^*$ we may assume
$r+q-b-1\geq 0$, namely $b\leq r+q-1$. By the assumption $q\leq n-r$ this
implies $b<n$. Hence Lem.(4-1)(3) completes the proof.
Next we show (1)$(ii)$. By what we have shown, we may assume $\ell\geq \e$.
By Lem.(4-1)(1) we may suppose $r+q-b-1\leq -s$, namely
$b\geq r+s+q-1$. The assumption $q\geq \frac{n-r+1}2$ implies
$r+s+q-1\geq (m-b)-(r+q-b-1)=n+s-q$. Hence Lem.(4-1)(2) completes
the proof.

To show that $\h q \ell$ is an isomorphism we show (a) and (b).
First assume $(2)(i)$. The assertion (b) has been shown in this
case. To show (a), by Lem.(4-1)$(1)^*$ we may assume
$r+q-a\geq 0$, namely $a\leq r+q$. By the assumption $q<n-r$ this
implies $a<n$. Hence Lem.(4-1)(3) completes the proof.
$(2)(ii)$ follows from $(2)(i)$ and Lem.(5-1)(1) by replacing $q$
by $n-r-q$ and $\ell$ by $\e-\ell$.
Next assume $(2)(iii)$. We only show (b). The proof of (a) is similar.
By Lem.(4-1)$(1),(2),(3)$ and $(1)^*,(2)^*,(3)^*$,
we have only  to consider either of the case $\ell=0$, $r+q-b-1=0$ and $b=n$
or the case $\ell=\e$, $r+q-b-1=-s$ and $b=r+s-1$. In the former case we have
$m-b=m-n=r+s-1<n$. Hence Lem.(4-1)(4) compltes the proof.
In the latter case we have $m-b=m-(r+s-1)=n\geq r+s$.
Hence Lem.(4-1)$(4)^*$ compltes the proof.
\par
Finally we show (3). By (1)$(ii)$ $\tau$ is injective.
By (1)$(i)$ and Lem.(5-2)(2), $\tau$ cannot be the zero map.
Thus $\tau$ is always an isomorphism.
\qed
\enddemo

\vskip 6pt
By Lem.(5-2) Theorem(II)(2) and (3) holds true in case $s=1$.
The case $s\geq 2$ is reduced to the special case by the induction on $s$
due to the following lemma(5-3). Let the notation be as \S 2. We put
$$\Sp=\Sigma_{\L}(-\log \sum^s_{j=2}\P_j) \qaq
\Sb=\Sigma_{\Lb}(-\log \sum^s_{j=2}\Pb_j),$$
where $\Lb=\L_{|\P_1}$ and $\Pb_j=\P_j\cap \P_1$ for $2\leq j\leq s$.
We also define the Jacobian rings
$$
\Rp q \l=\text{Coker}\big(
H^0(\Sp \ot \L^{q-1}\ot\pi^*\Cal O(\l)) \os{j(\sigma)}{\ra}
H^0(\L^{q}\ot\pi^*\Cal O(\l)) \big)
$$
$$
\Rb q \l=\text{Coker}\big(
H^0(\Sb \ot \Lb^{q-1}\ot\pi^*\Cal O(\l)) \os{j(\overline{\sigma})}{\ra}
H^0(\Lb^{q}\ot\pi^*\Cal O(\l)) \big)
$$
where
$\overline{\sigma}=\sum^r_{i=1} F_i \mu_i+\sum^s_{j=2} G_j \ld_j \in H^0(\Lb)$.
Put
$\d'=\d+e_1$ and $\e'=\e-e_1=\eb$.

\Lem 5-3. \it
(1) We have the exact sequence
$$\Rp {q-1} {\d'-n-1+\ell}@>\phi>> B_q(\d-n-1+\ell)
@>\rho>> \Rb q {\d-n-1+\ell}\to 0,$$
where $\rho$ is the reduction modulo $\ld_1\in A_q(-e_1)$ and $\phi$ is
the multiplication by $\ld_1$.
\vskip 4pt\noindent
(2) We have the exact sequence
$$\Rb {n-r-q} {\ell'-e_1}@>\psi>> B_{n-r-q}(\ell')
@>\pi>> \Rp {n-r-q} {\ell'}\to 0,$$
where $\pi$ is the natural projection arising from the natural injection
$\Sp\subset \Sigma$ and $\psi$ is the unique map which fits into the
commutative diagram
$$
\matrix
B_{n-r-q}(\ell'-e_1)&\\
\downarrow\rlap{$\rho$}&\searrow^{\tilde{\psi}}\\
\Rb {n-r-q}{\ell'-e_1}&@>\psi>>& B_{n-r-q}(\ell')\\
\endmatrix
$$
where $\tilde{\psi}$ is the multiplication by $G_1\in A_0(e_1)$.
\vskip 4pt\noindent
(3) The following diagram is commutative.
$$
\matrix
&&& 0\\
&&&\downarrow\\
&\Rp {q-1} {\d'-n-1+\ell}&@>{h'_{q-1}(\ell)}>>&
\Rp {n-(r+1)-(q-1)} {\d'+\e'-n-1-\ell}^*\\
&\downarrow\rlap{$\phi$} &&\downarrow\rlap{$\pi^*$}\\
&B_q(\d-n-1+\ell) &@>{h_q(\ell)}>> & B_{n-r-q}(\d+\e-n-1-\ell)^*\\
&\downarrow\rlap{$\rho$} &&\downarrow\rlap{$\psi^*$}\\
&\Rb q {\d-n-1+\ell}&@>{\overline{h}_q(\ell)}>>&
\Rb {n-r-q} {\d+\eb-n-1-\ell}^*\\
&\downarrow &&\\
&0 &&\\
\endmatrix$$
where the horizontal arrows are the duality maps defined before.

\demo{Proof}
The exactness of the sequences together with the well-definedness of $\psi$
is seen immediately from the explicit description of the Jacobian rings
(cf. Lem.(2-2)). The commutativity of the upper square of the diagram
in (3) is an easy consequence of the commutative diagram (cf. (5-2))
$$
\matrix
&\L^{-m+r+q-1}(\l)& \ra& \SSp \ot \L^{-m+r+q}(\l)&\ra \cdots \ra&
 \os{m}{\wedge}\SSp \ot \L^{r+q-1}(\l)&\ra&
\os{m+1}{\wedge}\SSp \ot \L^{r+q}(\l)\\
&\|&&\downarrow&&\downarrow&&\downarrow\\
 &\L^{-m+r+q-1}(\l)& \ra& \SS \ot \L^{-m+r+q}(\l)&\ra \cdots \ra&
 \os{m}{\wedge}\SS \ot \L^{r+q-1}(\l)&\ra&
\os{m+1}{\wedge}\SS \ot \L^{r+q}(\l)\\
\endmatrix
$$
where we put $\L(\ell)=\L\ot\pi^*\Cal O(\ell)$.
The vertical maps are the dual of the natural embedding
$\Sigma \hookrightarrow \Sp$.
\par
Next we show the commutativity of the lower sqaure.
By (5-1)(2) it suffices to show
$$\overline{\tau}(x)=\tau(\psi(x))\text{ for }x\in \Rb {n-r}{2(\d-n-1)+\eb},
\tag{5-3-1}$$
where $\taub$ and $\tau$ are the trace maps (cf. (5-3)).
We consider the following commutative diagram
$$
\matrix
&&& 0&& 0\\
&&&\uparrow&&\uparrow\\
&0&& \L^{n-r}(\l) &= & \L^{n-r}(\l) \\
&\uparrow&&\uparrow\rlap{$j'$}&&\uparrow\rlap{$j$}\\
0\gets&  \Lb^{n-r}(\l-e_1) &@<\alpha<<& \L^{n-r-1}(\l)\ot\Sp &@<\iota<< &
 \L^{n-r-1}(\l)\ot\Sigma &\gets 0\\
&\uparrow\rlap{$\overline{j}$}&&\uparrow&&\uparrow\\
0\gets& \Lb^{n-r-1}(\l-e_1)\ot\Sb &\gets&
\L^{n-r-1}(\l)\ot\os{2}{\wedge}\Sp &\gets &
 \L^{n-r-1}(\l)\ot\os{2}{\wedge}\Sigma &\gets 0\\
&\uparrow&&\uparrow&&\uparrow\\
& \vdots&& \vdots&& \vdots\\
&\uparrow&&\uparrow&&\uparrow\\
0\gets&  \Lb^{(-2r-s+1)}(\l-e_1)\ot\os{m}{\wedge}\Sb &\gets&
 \L^{-2r-s}(\l)\ot\os{m+1}{\wedge}\Sp &\gets &
 \L^{-2r-s}(\l)\ot\os{m+1}{\wedge}\Sigma &\gets 0\\
&\uparrow&&\uparrow&&\uparrow\\
&0&& 0&& 0\\
\endmatrix
$$
with $\ell=2(\d-n-1)+\e$. Here the vertical exact sequences come from (5-2)
and the horiontal exact sequences come from the exact sequence
$$ 0\to\Sigma \to \Sp \to N_{\P_1/\P}\to 0$$
(coming from the exact sequence
$0\to T_{\P}(-\log \sum^s_{j=1}\Pb_j)\to
 T_{\P}(-\log \sum^s_{j=2}\Pb_j)\to N_{\P_1/\P}\to 0$) and the isomorphism
$$N_{\P_1/\P}\simeq \Cal O_{\P}(\P)\ot\Cal O_{\P_1} \simeq
\L(-e_1)\ot\Cal O_{\P_1}.$$
We note that the bottom row is isomorphic to the adjunction sequence
$$ 0\gets K_{\P_1} \gets K_{\P}\ot\Cal O_{\P}(\P_1)\gets K_{\P}\gets 0$$
and the left and right vetical sequences induce the maps
$$\text{Coker}\big(
H^0(\L^{n-r-1}(\ell)\ot\Sigma) @>{j}>>H^0(\L^{n-r}(\ell))\big)\to
H^{m}(\P,K_{\P})\simeq k,$$
$$\text{Coker}\big(
H^0(\Lb^{n-r-1}(\ell-e_1)\ot\Sb) @>\overline{j}>>H^0(\Lb^{n-r}(\ell-e_1))
\big)\to H^{m-1}(\P_1,K_{\P_1})\simeq k,$$
which are nothing but the trace maps $\tau$ and $\taub$ respectively.
On the other hand we note that the map
$$ \alpha\sp:\sp H^0(\L^{n-r-1}(\ell)\ot\Sp) \to H^0(\Lb^{n-r}(\ell-e_1))$$
is surjective. Using this we define the map
$$ \delta\sp:\sp
H^0(\Lb^{n-r}(\l-e_1))\to
\text{Coker}
\big(H^0(\L^{n-r-1}(\ell)\ot\Sigma) @>{j}>>H^0(\L^{n-r}(\ell))\big)
$$
by $\delta(x)=j'(\tilde x) \text{ mod Im}(j)$ for
$x\in H^0(\Lb^{n-r}(\ell-e_1))$ and
$\tilde x\in H^0(\L^{n-r-1}(\ell)\ot\Sp)$ with $\alpha(\tilde{x})=x$.
It is easily seen that $\delta$ coincides with the multiplication by
$G_1\in H^0(\Cal O_{\P^n}(e_1))$.
Thus, to show (5-3-1) it suffices to prove
\def\part{\partial}
$\part_1(x)=\part_2(j'(\tilde x))$, where
$$ \part_1\sp:\sp H^0(\Lb^{n-r}(\ell-e_1))\to H^1(\Ker(j))\qaq
 \part_2\sp:\sp H^0(\L^{n-r}(\ell))\to H^1(\Ker(j))$$
are the boundary maps coming from the exact sequences
$$0\to\Ker(j)@>\iota>>\Ker(j')@>\alpha>> \Lb^{n-r}(\ell-e_1)\to 0,\quad
0\to\Ker(j)\to \L^{n-r-1}(\ell)\ot \Sigma @>j>> \L^{n-r}(\ell)\to 0.$$
Take
$$\text{an open covering } \P=\cup_{i\in I} U_i\qaq
 \{\eta_i\}_{i\in I}\in \prod_{i\in I} H^0(U_i,\L^{n-r-1}(\ell)\ot \Sigma)$$
such that $j(\eta_i)=j'(\tilde{x})_{|U_i}$. Putting
$\xi_i=\tilde{x}_{|U_i} -\iota(\eta_i)$, we see
$ \eta_i\in H^0(U_i,\Ker(j'))$ and $\alpha(\xi_i)= x_{|U_i}$. Thus
$$ \eta_{i|U_i\cap U_j} - \eta_{j|U_i\cap U_j}=
 \xi_{i|U_i\cap U_j} - \xi_{j|U_i\cap U_j}\in
H^0(U_i\cap U_j,\Ker(j))$$
is a Cech cocycle representing both $\part_1(x)$ and $\part_2(j'(\tilde{x}))$.
This completes the proof.
\qed
\enddemo

\vskip 6pt
Finally we prove Theorem(II)(2)$(iii)$.
We reduce it to the case $s=1$. For this we consider
the diagram of Lem.(5-3) in case $s=1$ and $\ell=0$.
$$
\matrix
&&& 0\\
&&&\downarrow\\
&\Rp {q-1} {\d'-n-1}&@>{h'}>>&
\Rp {n-(r+1)-(q-1)} {\d'-n-1}^*\\
&\downarrow\rlap{$\phi$} &&\downarrow\rlap{$\pi^*$}\\
&B_q(\d-n-1) &@>{h}>> & B_{n-r-q}(\d+e_1-n-1)^*\\
&\downarrow\rlap{$\rho$} &&\downarrow\rlap{$\psi^*$}\\
&\Rb q {\d-n-1}&@>{\overline{h}}>>&
\Rb {n-r-q} {\d-n-1}^*\\
&\downarrow &&\\
&0 &&\\
\endmatrix$$
We want to show $\overline{h}$ is an isomorphism.
By Theorem(I) $B_q(\d-n-1)$ and $B_{n-r-q}(\d-n-1)$ have the same dimension.
Hence it suffices to show
the injctivity of $\overline{h}$.
First we have the following.

\vskip 6pt \medbreak\noindent
\it Claim. \it $\Rb q {\d-n-1}=0$ if $q>n-r\geq 1$.
\rm\demo{Proof}
By Lem.(5-2)(2)$(iii)$, $B_q(\d-n-1)=0$ if $q>n-r\geq 1$.
Hence the claim follows from the surjectivity of $\rho$.
\qed
\enddemo
\vskip 5pt

Assume $n-r\geq 1$.
By Lem.(5-2) $h$ is an isomorphism so that $\overline{h}$
is injective if $h'$ is surjective. By the claim and by
the induction on $r$ we are reduced to show the injectivity of $\overline{h}$
in case $n-r=1$ and $q=1$ in which case $\psi^*$ is surjective. For its dual
$\psi:\Rb 0 {\d-n-1}\to B_0(\d+e_1-n-1)$ that is the multiplication by $G_1$,
is injective since $(F_1,\dots, F_r,G_1)$ is a regular sequence in
$k[X_0,\dots,X_n]$. Since $h$ is an isomorphism by Lem.(5-2), the diagram
shows that $\overline{h}$ is surjective so that injective by the reason of
dimension. This completes the proof of Theorem(II)(2)$(iii)$ in case
$n-r\geq 1$. The diagram implies further that $h'$ is
an isomorphism in case $n-r=1$ and $q=1$ so that Theorem(II)(2)$(iii)$
in case $n-r=0$ is also proved.
\qed


\vskip 20pt

\input amstex
\documentstyle{amsppt}
\hsize=16cm
\vsize=23cm

\head \S6. Proof of Theorem(II'). \endhead
\vskip 8pt

In this section we prove Theorem(II').
First the surjectivity of $\eXZ$ follows from Thoerem(II)(3).
As for $\Ker(\eXZ)$ we first show the following.

\Pr 6-1. \it
$\wedge_X^{n-r}(G_1,\dots,G_s)\subset \Ker(\eXZ).$
\rm\vskip 6pt

Note that $B_0(\d+\e-n-1)=\Pol^{\d+\e-n-1}/(F_1,\dots,F_r)$
where $\Pol^\ell\subset k[X_0,\dots,X_n]$ is the subspace of homogeneous
polynomials of degree $\ell$ and $(F_1,\dots,F_r)\subset \Pol^\ell$ is
the subspace generated by the multiples of $F_i$.
By Thoerem(I) we have the isomorphisms
$$ B_0(\d+\e-n-1) \isom H^0(X,\WXZ {n-r});
\spa A\to \Res X \frac{A}{F_1\cdots F_r G_1\cdots G_s} \Omega,$$
where
$$ \Omega:= \sum_{i=0}^n (-1)^i X_i dX_0\wedge\cdots \wedge\widehat{dX_i}\wedge
 \cdots\wedge dX_n,$$
and
$$ \Res X\scs \Omega_{\P^n}^n(\log X_*+Y_*) \to \WXZ {n-r}$$
is the composite of the residue maps along $F_1=F_2=\cdots=F_r=0$.

\Def 6-1. \it Assume $s\geq n-r+1$.
For integers $1\leq j_1<\cdots< j_{n-r+1}\leq s$, we write
$$ A(j_1,\dots,j_{n-r+1}):= \sum_{\sigma\in \pg {n+1}}\sign(\sigma)
\frac{\partial F_1}{\partial X_{\sigma(0)}}\cdots
\frac{\partial F_r}{\partial X_{\sigma(r-1)}}\cdot
\frac{\partial G_{j_1}}{\partial X_{\sigma(r)}}\cdots
\frac{\partial G_{j_{n-r+1}}}{\partial X_{\sigma(n)}},$$
$$ A'(j_1,\dots,j_{n-r+1}):= A(j_1,\dots,j_{n-r+1})\cdot
\frac{G_1\cdots G_s}{G_{j_1}\cdots G_{j_{n-r+1}}}\in S^{\d+\e-n-1}.$$
\rm\vskip 4pt

\Lem 6-1. \it We have (cf. Def.(1-3))
$$ \Res X \frac{A'(j_1,\dots,j_{n-r+1})}{F_1\cdots F_r G_1\cdots G_s} \Omega
=(-1)^{r+1}\omega_X(j_1,\dots,j_{n-r+1}).$$
\rm

\Lem 6-2. \it Write $A'=A'(j_1,\dots,j_{n-r+1})$. Then we have (cf. Def.(1-1))
$$A'\mu_i, \spa A'\lambda_j \in J(\ul{F},\ul{G}) \qfor
1\leq \forall  i\leq r, \spa 1\leq \forall  j\leq s.$$

\rm\vskip 5pt

Pr.(6-1) follows from the above lemmas:
By Lem.(6-1) it suffices to show
$$A'(j_1,\dots,j_{n-r+1})\in
\Ker( B_0(\d+\e-n-1) @>{{\h {n-r} 0}^*}>> B_{n-r}(\d-n-1)^*).$$
Since ${\h {n-r} 0}^*$ is given by the pairing
$$B_0(\d+\e-n-1) \otimes B_{n-r}(\d-n-1)\to B_{n-r}(2(\d-n-1)+\e) @>\tau>> k,$$
Pr.(6-1) follows from Lem.(6-2).

\vskip 5pt\noindent
\it {Proof of Lem.(6-1)} \rm
We may suppose $j_1=1,\dots,j_{n-r+1}=n-r+1$. We may prove the formula
on the affine subspace $\{X_0\not=0\}$. Let
$$ A= \sum_{\sigma\in \pg {n+1}}\sign(\sigma)
\frac{\partial F_1}{\partial X_{\sigma(0)}}\cdots
\frac{\partial F_r}{\partial X_{\sigma(r-1)}}\cdot
\frac{\partial G_{1}}{\partial X_{\sigma(r)}}\cdots
\frac{\partial G_{{n-r+1}}}{\partial X_{\sigma(n)}}.$$
For polynomials $h_1,\dots, h_n\in k[X_0,\dots, X_n]$ write
$$ J(h_1,\dots, h_n)=
\det\pmatrix
\frac{\partial h_1}{\partial X_1} &\hdots & \frac{\partial h_1}{\partial X_n}\\
\vdots&&\vdots\\
\frac{\partial h_n}{\partial X_1} &\hdots & \frac{\partial h_n}{\partial X_n}\\
\endpmatrix.$$
Writing $G_j=F_{r+j}$ and $e_j=d_{r+j}$ for $1\leq j\leq s$, we claim
$$ X_0 A=
\sum_{\nu=1}^n (-1)^{\nu-1} (d_\nu\cdot F_\nu)
J(F_1,\dots,\widehat{F_\nu},\dots,F_{n+1})\leqno(*)$$
that implies that we have on $\{X_0\not=0\}$
$$\multline
\Res X \frac{A'(j_1,\dots,j_{n-r+1})}{F_1\cdots F_r G_1\cdots G_s} \Omega=
\Res X \sum_{\nu=1}^r (-1)^{\nu+1} d_\nu
\frac{df_1}{f_1}\wedge\cdots\wedge
\widehat{\frac{df_\nu}{f_\nu}}
\wedge\cdots\wedge\frac{df_r}{f_r}\wedge
\frac{dg_1}{g_1}\wedge\cdots\wedge
\frac{dg_{{n-r+1}}}{g_{{n-r+1}}}\\
+ \Res X \sum_{\mu=1}^{n-r+1}(-1)^{r+\mu+1} e_{{\mu}}
\frac{df_1}{f_1}\wedge\cdots\wedge\frac{df_r}{f_r}\wedge
\frac{dg_1}{g_1}\wedge\cdots\wedge\widehat{\frac{dg_\mu}{g_\mu}}
\wedge\cdots\wedge\frac{dg_{n-r+1}}{g_{n-r+1}}.
\endmultline$$
where $f_i=F_i/X_0^{d_i}$ and $g_j=G_j/X_0^{e_j}$.
Since the first term vanishes this completes the proof of Lem.(6-1).
To show the formula $(*)$ we note
$$ X_0\frac{\partial F_\nu}{\partial X_0}=
d_\nu\cdot F_\nu -\sum_{i=1}^n X_i\frac{\partial F_\nu}{\partial X_i}.$$
We have
$$\eqalign{
X_0 A&=\sum_{\sigma\in \pg {n+1}}\sign(\sigma)
X_0\frac{\partial F_1}{\partial X_\sigma(0)}\cdots
\frac{\partial F_{n+1}}{\partial X_\sigma(n)}\cr
&=\sum_{\sigma(0)=0}(*)+\sum_{\sigma(1)=0}(*)+\cdots+\sum_{\sigma(n)=0}(*)\cr
&=\sum_{\sigma\in \pg {n}}\sign(\sigma)
(d_1\cdot F_1-\sum_{i=1}^n X_i\frac{\partial F_1}{\partial X_i})\cdot
\frac{\partial F_{2}}{\partial X_{\sigma(1)}}
\cdots
\frac{\partial F_{n+1}}{\partial X_{\sigma(n)}}\cr
&-\sum_{\sigma\in \pg {n}}\sign(\sigma)
(d_2\cdot F_2-\sum_{i=1}^n X_i\frac{\partial F_2}{\partial X_i})\cdot
\frac{\partial F_{1}}{\partial X_{\sigma(1)}}
\frac{\partial F_{3}}{\partial X_{\sigma(2)}}
\cdots
\frac{\partial F_{n+1}}{\partial X_{\sigma(n)}}+\cdots\cr
&+(-1)^{n+1}\sum_{\sigma\in \pg {n}}\sign(\sigma)
(d_{n+1}\cdot F_{n+1}-\sum_{i=1}^n X_i\frac{\partial F_2}{\partial X_i})\cdot
\frac{\partial F_{1}}{\partial X_{\sigma(1)}}
\cdots
\frac{\partial F_{n}}{\partial X_{\sigma(n)}}\cr
&=\sum_{\nu=1}^n (-1)^{\nu-1} (d_\nu\cdot F_\nu)
J(F_1,\dots,\widehat{F_\nu},\dots,F_{n+1})
-\sum_{i=1}^n X_i P_i\cr}$$
where
$$P_i=\sum_{\nu=1}^{n+1} (-1)^{\nu+1}
\frac{\partial F_\nu}{\partial X_i}\cdot
 J(F_1,\dots,\widehat{F_\nu},\dots,F_{n+1})
=\det\pmatrix
\frac{\partial F_1}{\partial X_i}&
\frac{\partial F_1}{\partial X_1}&\hdots &\frac{\partial F_1}{\partial X_n}\\
\frac{\partial F_2}{\partial X_i}&
\frac{\partial F_2}{\partial X_1}&\hdots &\frac{\partial F_2}{\partial X_n}\\
\vdots&\vdots&&\vdots&\\
\frac{\partial F_{n+1}}{\partial X_i}&
\frac{\partial F_{n+1}}{\partial X_1}&\hdots
&\frac{\partial F_{n+1}}{\partial X_n}\\
\endpmatrix =0.$$
This completes the proof.
\qed

\vskip 6pt\noindent
\it {Proof of Lem.(6-2)} \rm
We may suppose $j_1=1,\dots,j_{n-r+1}=n-r+1$.
Modulo $J(\ul{F},\ul{G})$ we have
$$ \eqalign{A'\mu_1&= \sum_{\sigma\in \pg {n+1}}\sign(\sigma)
\frac{\partial F_1}{\partial X_{\sigma(0)}}\mu_1\cdot
\frac{\partial F_2}{\partial X_{\sigma(1)}}\cdots
\frac{\partial F_r}{\partial X_{\sigma(r-1)}}\cdot
\frac{\partial G_1}{\partial X_{\sigma(r)}}\cdots
\frac{\partial G_{n-r+1}}{\partial X_{\sigma(n)}}
\cdot G_{n-r+2}\cdots G_s\cr
&\equiv -\sum_{\sigma\in \pg {n+1}}\sign(\sigma)
(\frac{\partial F_2}{\partial X_{\sigma(0)}}\mu_2+\cdots+
\frac{\partial G_s}{\partial X_{\sigma(0)}}\lambda_s)\cdot
\frac{\partial F_2}{\partial X_{\sigma(1)}}\cdots
\frac{\partial G_{n-r+1}}{\partial X_{\sigma(n)}}
\cdot G_{n-r+2}\cdots G_s \cr
&\equiv -\sum_{\sigma\in \pg {n+1}}\sign(\sigma)
(\frac{\partial F_2}{\partial X_{\sigma(0)}}\mu_2+\cdots+
\frac{\partial G_{n-r+1}}{\partial X_{\sigma(0)}}\lambda_{n-r+1})\cdot
\frac{\partial F_2}{\partial X_{\sigma(1)}}\cdots
\frac{\partial G_{n-r+1}}{\partial X_{\sigma(n)}}
\cdot G_{n-r+2}\cdots G_s \cr}$$
The coefficient of $\mu_i$ $(2\leq i\leq r)$ in the above is
$$
-(G_{n-r+2}\cdots G_s)\sum_{\sigma\in \pg {n+1}}\sign(\sigma)
\frac{\partial F_i}{\partial X_{\sigma(0)}}\cdot
\frac{\partial F_2}{\partial X_{\sigma(1)}}\cdots
\frac{\partial G_{n-r+1}}{\partial X_{\sigma(n)}}
=0$$
Similarly the coefficient of $\lambda_j$ $(1\leq j\leq n-r+1)$ vanishes.
This proves $A'\mu_1\equiv 0\mod J(\ul{F},\ul{G})$.
The rest of the assertion is proven in the same manner.
\qed
\vskip 6pt

Due to Pr.(6-1) Th.(II') now follows from
$$ \dim_k(\Ker(\eXZ))=\dim_k({\h {n-r} 0}^*)
\leq \dim_k\wedge_X^{n-r}(G_1,\dots,G_s).\leqno(*)$$
Note that by Th.(II)(2), $\dim_k({\h {n-r} 0}^*)=0$ if $s\leq n-r$ and
$n-r\geq 1$.
Consider the following commutative diagram that is the dual of the diagram
in Lem.(5-3)
$$\matrix
&& 0\\
&&\downarrow\\
\Rb 0 {\d+\eb-n-1} &@>{\overline{h}_{n-r}(0)^*}>> & \Rb {n-r} {\d-n-1}^*\\
\downarrow\rlap{$\psi$} &&\downarrow\rlap{$\rho^*$}\\
B_0(\d+\e-n-1) &@>{h_{n-r}(0)^*}>> & B_{n-r}(\d-n-1)^*\\
\downarrow\rlap{$\pi$} &&\downarrow\rlap{$\phi^*$}\\
\Rp 0 {\d'+\e'-n-1} &@>{h'_{n-(r+1)}(0)^*}>> & \Rp {n-(r+1)} {\d'-n-1}^*\\
\downarrow\\
0 \\
\endmatrix\leqno(6-1)$$
where the horizontal maps are surjective by Theorem(II)(3).
Combining this with the exact sequence
$$  0\to \wedge_X^{n-r}(G_2,\dots,G_{s})\to \wedge_X^{n-r}(G_1,\dots,G_s)
@>{\Res {Z_1}}>> \wedge_{Z_1}^{n-r-1}(G_2,\dots,G_s)\to 0,$$
the assertion $(*)$ is reduced by induction on $s$ to the case $n-r=1$.
Since $\dim_k\wedge_X^{1}(G_1,\dots,G_s)=s-1$ as is easily seen,
it follows by induction on $s$ from (6-1) and the following.

\Lem 6-3. \it Assuming $n=r$ and $s\geq 1$, $\dim_k\Ker(h_0(0)^*)=1$ where
$$ h_0(0)^*\scs B_0(\d+\e-n-1) \to B_0(\d-n-1)^*.$$
\rm
\demo{Proof}
We have the following commutative diagram
$$\matrix
0\\
\downarrow\\
S^{\d+\eb-n-1}/(F_1,\dots,F_n) &=& \Rb 0 {\d+\eb-n-1}
&@>{\overline{h}_{n-r}(0)^*}>> & \Rb {0} {\d-n-1}^*\\
\downarrow\rlap{$G_1$}&&\downarrow\rlap{$\psi$} &&\|\\
S^{\d+\e-n-1}/(F_1,\dots,F_n) &=& B_0(\d+\e-n-1)
&@>{h_{n-r}(0)^*}>> & B_{0}(\d-n-1)^*\\
\downarrow\\
S^{\d+\e-n-1}/(F_1,\dots,F_n,G_1)\\
\downarrow\\
0\\
\endmatrix\leqno(6-2)$$
The left vertical sequence is exact due to the assumption (1-1) in \S1.
The induction hypothesis and Theorem(II)(2)$(ii)$ imply
$$ \dim_k(\overline{h}_{n-r}(0)^*)=\left\{\aligned
& 0\quad\text{ if } s=1,\\
& 1\quad\text{ if } s\geq 2.\\
\endaligned\right.$$
Thus the lemma follows by induction on $s$ from (6-2) and the following fact
that is a consequence of Macaulay's theorem (cf. [Do, Th.2.5])
$$ \dim_k S^{\d+\e-n-1}/(F_1,\dots,F_n,G_1)=
\left.\left\{\gathered
 1\\
  0
\endgathered\right.\qquad
\aligned
&\text{if $s=1$,}\\
&\text{if $s\geq 2$}
\endaligned\right.
$$
\qed
\enddemo


\vskip 20pt

\input amstex
\documentstyle{amsppt}
\hsize=16cm
\vsize=23cm

\head \S7. Proof of Theorem(III) \endhead
\vskip 8pt

In this section we complete the proof of Th.(III).We deduce it from the
following.

\Th 7-1. \it Assume $s\geq 1$.
Let $W \subset A_1(0)$ is a base point free subspace of codimension $c$
({\it i.e.} for any $x \in \P^n(\Bbb C)$, the evaluation map
$W \subset A_1(0) \ra \us{i}{\op} \Bbb C \mu_i \op \us{j}{\op}\Bbb C \lambda_j$ at
$x$ is surjective). Then the Koszul complex
$$
B_p(\l) \otimes \os{q+1}{\wedge}W \rightarrow
B_{p+1}(\l) \otimes \os{q}{\wedge}W \rightarrow
B_{p+2}(\l) \otimes \os{q-1}{\wedge}W
$$
is exact if one of the following conditions is  satisfied.
\roster
\item"$(i)$"
$p\geq 0$, $q=0$ and $\md p+\ell\geq c$.
\item"$(ii)$"
$p\geq 0$, $q=1$ and $\md p+\ell\geq 1+c$ and
$\md(p+1)+\ell\geq d_{max}+c$.
\item"$(iii)$"
$p\geq 0$, $\md(r+p)+\ell-\d\geq q+c$, $\ell\geq \d-n-1$,
$e_1=\cdots=e_s$ and either $r+s\leq n+2$ or $p\not=n-r-1$.
\item"$(iv)$"
$p\geq 0$, $\md(r+p)+\ell-\d\geq q+c$, $\d+e_{max}-n-1>\ell\geq \d-n-1$
and either $r+s\leq n+2$ or $p\not=n-r-1$.
\endroster
\rm
\vskip 6pt

First we deduce Th.(III) from Th.(7-1).
Let $W:=\Ker(A_1(0) \ra B_1(0)/V)$.
Since $W$ contains $J:=J(X,\Zst) \cap A_1(0)$ (cf. Def.(1-2)),
it is a base point free subspace of codimension $c$.
We have the Koszul exact sequence
$$
0 \rightarrow
S^{\cdot}(J)
\rightarrow
W \ot S^{\cdot-1}(J)
\rightarrow
\cdots \rightarrow
\os{\cdot-1}{\wedge} W \ot J
\rightarrow
\os{\cdot}{\wedge} W
\rightarrow
\os{\cdot}{\wedge} V
\rightarrow 0.
$$
This complex tensored with $B_*(\l)$ induces the following diagram.
$$
\matrix
\cdots                                   & \rightarrow &
B_p(\l)\ot\os{q+1-i}{\wedge}W\ot S^i(J)&
\rightarrow                              & \cdots      & \rightarrow &
B_p(\l)\ot\os{q+1}{\wedge}W            & \rightarrow
B_p(\l)\ot\os{q+1}{\wedge}V    &
\rightarrow                              & 0 \\
                                         &             &
\downarrow                               &
                                         &             &             &
\downarrow                               &
\downarrow                               &
                                         & \\
\cdots                                   & \rightarrow &
B_{p+1}(\l)\ot\os{q-i}{\wedge}W\ot S^i(J)&
\rightarrow                              & \cdots      & \rightarrow &
B_{p+1}(\l)\ot\os{q}{\wedge}W            & \rightarrow
B_{p+1}(\l)\ot\os{q}{\wedge}V    &
\rightarrow                              & 0 \\
                                         &             &
\downarrow                               &
                                         &             &             &
\downarrow                               &
\downarrow                               &
                                         & \\
\cdots                                   & \rightarrow &
B_{p+2}(\l)\ot\os{q-1-i}{\wedge}W\ot S^i(J)&
\rightarrow                              & \cdots      & \rightarrow &
B_{p+2}(\l)\ot\os{q-1}{\wedge}W            & \rightarrow
B_{p+2}(\l)\ot\os{q-1}{\wedge}V    &
\rightarrow                              & 0 \\
                                         &             &
\downarrow                               &
                                         &             &             &
\downarrow                               &
\downarrow                               &
                                         & \\
                                         &             &
\vdots                                   &
                                         &             &             &
\vdots                                   &
\vdots                                   &
                                         &
\endmatrix
$$
where the vertical sequences are the Koszul complexes tensored with
$S^{\cdot}(J)$. Since $J$ annihilates $B_p(\ell)$, the diagram is
commutative. Therefore to show the exactness of the complex in Th.(III),
it suffices to show that
$$
B_{p+i}(\l) \ot \os{q-2i+1}{\wedge} W
\lra
B_{p+i+1} (\l)\ot \os{q-2i}{\wedge} W
\lra
B_{p+i+2}(\l) \ot \os{q-2i-1}{\wedge} W
$$
is exact for $\forall i\geq 0$ and it follows from Th.(7-1) under the
assumptions of Th.(III).
\qed

\vskip 8pt

For the proof of Th.(7-1) we recall the regularity of sheaves ([G2]).
A coherent sheaf $\Cal F$ on $\Bbb P^n$ is called $m$-{\it regular} if
$$
H^i(\Bbb P^n,\Cal F \ot \Cal O_{\Bbb P^n}(m-i))=0 \quad
             \text{for $\forall i>0$.}
$$
We use the following properties of the regularity of sheaves, whose proof
can be found in [G2].
\roster
 \item If $\Cal F$ is $m$-regular, then also $(m+1)$-regular.
 \item If $\Cal F$ and $\Cal F'$ are $m$-regular and $m'$-regular respectively, then $\Cal F \ot \Cal F'$ is $(m+m')$-regular.
\endroster
In particular, if $E$ is a $m$-regular locally free sheaf on $\Bbb P^n$,
then $\os{p}{\wedge} E$ is $(mp)$-regular since it is a direct summand of
$E^{\ot p}$.
Let $\l \geq 0$ be an integer, and define a locally free sheaf $E$ on a projective space $\P^n$ by the exact sequence
$$
0 \lra E \lra
H^0(\P^n,\cO_{\P^n} (\l)) \otk \Cal O_{\P^n}\lra
\Cal O_{\Bbb P^n}(\l) \lra 0.
$$
Then clearly $E$ is 1-regular, therefore $\os{p}{\wedge}E$ is $p$-regular.
In [G2], there is a further result: We replace $H^0(\Cal O_{\Bbb P^n}(\l))$
by $V$ a base point free linear subspace of $H^0(\Cal O_{\Bbb P^n}(\l))$ of
codimension $c$ and define $E'$ by
$$
0 \lra E' \lra
V \otk \Cal O_{\P^n} \lra \Cal O_{\Bbb P^n}(\l) \lra 0,
$$
then $\os{p}{\wedge}E'$ is $(p+c)$-regular.
This argument is applicable not only to $\Cal O_{\Bbb P^n}(\l)$ but also to
any locally free sheaf satisfying certain conditions. We will need it later.

\Lem7-1. \it
Let $\Cal N$ be a locally free sheaf on $\P^n$ generated by global sections. We assume that $\Cal N$ satisfies $H^p(\Cal N(-p))=0$ for $0<p<n$ ({\it e.g.} $\Cal N=\E$).
Let $V$ be a linear subspace of $H^0(\Cal N)$ of codimension $c$, such that
$V \otk \Cal O_{\P^n} \ra \Cal N$ is surjective ({\it i.e.} base point free). Define the locally free sheaf $N$ by the exact sequence
$$
0 \lra N \lra
V \otk \Cal O_{\P^n} \lra \Cal N \lra 0.
$$
Then $\os{p}{\wedge}N$ is $(p+c)$-regular. \rm
\vskip 6pt

Now we go back to the proof of Th.(7-1). The following lemma is a
generalization of [G2, Th.4.1].

\Lem7-2. \it
Let $q \geq 0$, $\nu \geq 0$, $\ell$ integers. Then the Koszul complex
$$
A_{\nu}(\l) \ot \os{q+1}{\wedge}W  \lra
A_{\nu+1}(\l) \ot \os{q}{\wedge}W  \lra
A_{\nu+2}(\l) \ot \os{q-1}{\wedge}W
$$
is exact if $\md\nu +\l \geq c+q$. \rm

\demo{Proof}
We define a locally free sheaf $M$ on $\P$ by the exact sequence
$$
0 \lra M \lra
W \otk \Cal O_{\P} \lra \Cal L \lra 0.
\tag *
$$
where the first map comes from the identification
$H^0(\P,\cL)=A_1(0)\supset W$ (cf. Lem.(2-2)). Then we obtain a Koszul exact
sequence
$$
0 \ra \os{q+1}{\wedge} M
  \ra \os{q+1}{\wedge} W \otk \Cal O_{\Bbb P}
  \ra \os{q}{\wedge} W \otk \Cal L
  \ra \cdots
  \ra W\otk \cL^q
  \ra \Cal L^{q+1}
  \ra 0.
$$
Tensoring with
$\Cal L^{\nu} \ot \pi^* \Cal O_{\Bbb P^n}(\ell)$,
this gives an acyclic resolution of
$
\os{q+1}{\wedge} M \ot \Cal L^{\nu} \ot \pi^* \Cal O_{\Bbb P^n}(\ell)
$.
By Lem.(2-1)
$$
H^i(\Bbb P, \Cal L^{\nu^{\prime}} \ot \pi^* \Cal O_{\bP^n}(\ell))
\simeq
H^i(\Bbb P^n, S^{\nu^{\prime}}(\Cal E) \ot \Cal O_{\Bbb P^n}(\l))
=\underset{\alpha}\to{\op} H^j(\Bbb P^n, \Cal O_{\Bbb P^n}(\alpha))
\qfor \forall \nu'\geq \nu
$$
with $\alpha\geq \md\nu^{\prime}+\l \geq \md\nu+\l \geq c+q \geq 0$ so that
it vanishes if $i>0$.
Therefore the cohomology group in Lem.(7-2) is isomorphic to
$$H^1(\os{q+1}{\wedge} M \ot \Cal L^{\nu} \ot \pi^* \Cal O_{\Bbb P^n}(\ell))$$
which we shall prove vanishes. Let $\pi:\bP \to \bP^n$ be the projection.
We apply $\pi_*$ to $(*)$ and get the exact sequence
$$
0 \lra \pi_*M \lra
W' \otk \Cal O_{\P^n} \lra \Cal E \lra 0
\qwith W'=\pi_* W\subset H^0(\bP^n,\cE).
$$
The surjectivity of the right map is due to the base point freeness of $W$.
Put $N=\pi_*M$. Then by Lem.(7-1), $\os{i}{\wedge}N$ is $(c+i)$-regular.
On the other hand we have the commutative diagram:
$$
\matrix
0 \lra & \pi^* N
 &\lra & W \otk \Cal O_{\Bbb P}
 &\lra & \pi^*\Cal E
 &\lra 0 \\
       & \quad \downarrow g
 &     & \quad \downarrow =
 &     & \quad \downarrow g'
 &       \\
0 \lra & M
 &\lra & W \otk \Cal O_{\Bbb P}
 &\lra & \Cal L
 &\lra 0.
\endmatrix
$$
where the vertical maps are induced by the adjunction for $\pi$.
By the snake lemma, $g$ is injective and $\Coker(g) \simeq\Ker(g')$.
By the exact sequence (cf. Lem.(2-1)(3))
$$ 0\to \Omega_{\P/\P^n} \to \pi^*\cE\otimes \cL^{-1}\to \cO_{\P}\to 0,
\leqno(**)$$
we get $\Ker(g') \simeq \Omega^1_{\Bbb P/\Bbb P^n} \ot \Cal L$.
Hence we have the exact sequence
$$
0 \lra  \pi^* N
  \lra M
  \lra \Omega^1_{\Bbb P/\Bbb P^n} \ot \Cal L
  \lra 0.
$$
which induces the filtration
$$
\os{q+1}{\wedge} M
  = F^0 \supset F^1 \supset \cdots \supset F^{q+1} \supset F^{q+2}=0
$$
such that
$\text{Gr}^i_F(\os{q+1}{\wedge} M)=F^i/F^{i+1}
   \simeq \pi^*(
       \os{i}{\wedge}N
   )\ot \Omega^{q-i+1}_{\Bbb P/\Bbb P^n} \ot \Cal L^{q-i+1}$.
So it suffices to show that
$$
H^1(\Bbb P, \Cal L^{\nu+q-i+1} \ot \Omega^{q-i+1}_{\Bbb P/\Bbb P^n}
            \ot \pi^*(\os{i}{\wedge}N
            \ot \Cal O_{\Bbb P^n}(\l)))=0\quad
\text{for } 0 \leq \forall i \leq q+1
$$
The exact sequence $(**)$ induces the exact sequence
$$
0 \ra \os{r+s}{\wedge}\pi^*\E \ot \L^{-r-s}
\ra \cdots
\ra  \os{p+1}{\wedge}\pi^*\E \ot \L^{-p-1}
\ra  \Om_{\P/\P^n}^p
\ra 0.
$$
Therefore it suffices to show that
$$
H^j(\Bbb P, \Cal L^{\nu-j} \ot \pi^*(
            \os{q+1-i+j}{\wedge} \Cal E
            \ot \os{i}{\wedge}N
            \ot \Cal O_{\Bbb P^n}(\l)))=0\quad
\text{for }1 \leq \forall j \leq r+s-(q+1-i)\text{ and }0\leq\forall i\leq q+1.
$$

In case $1 \leq j \leq \nu $ the above cohomology is isomorphic to
$$
H^j(\Bbb P^n, S^{\nu-j}(\Cal E)
            \ot \os{q+1-i+j}{\wedge} \Cal E
            \ot \os{i}{\wedge}N
            \ot \Cal O_{\Bbb P^n}(\l))
\simeq
\underset{\alpha}\to{\op} H^j(\Bbb P^n, \Cal O_{\Bbb P^n}(\alpha)
            \ot \os{i}{\wedge}N)
$$
with
$\alpha \geq \md(\nu-j)+\md(q+1-i+j)+\l =\md(\nu+q+1-i)+\l$.
Since $\os{i}{\wedge}N$ is $(i+c)$-regular, this vanishes
if $\alpha+j \geq i+c$, which holds since
$$\alpha+j-i-c \geq \md(\nu+q+1-i)+\l+j-i-c\geq \md\nu+\l-(q+c) \geq 0.$$
\par
Next assume $j >\nu $.
If $\nu-j>-r-s$, the cohomology vanishes by Lem.(2-1)(2).
Hence we only consider the case $\nu-j \leq -r-s$, namely $j \geq r+s+\nu$.
Since $j \leq r+s-(q+1-i)$ by the assumption, we only consider the case
$\nu=0$, $j=r+s$, $i=q+1$. Then the cohomology is isomorphic to
$H^1(\Bbb P^n, \Cal O_{\Bbb P^n}(\l) \ot \os{q+1}{\wedge}N)$.
Since $\os{q+1}{\wedge}N$ is ($q+1+c$)-regular and
$\l+1= \md\nu+\l+1 \geq q+1+c$, it vanishes. This completes the proof of
Lem.(7-2).
\qed
\enddemo

Now we prove Th.(7-1).
Write $\Sigma=\Sigma_{\cL}(\log \P_*)$ and put
$$\cM_{k,h}(\l)=\os{m+1-h}{\wedge}\SS \ot \L^{r+k-h}
\ot \pi^*\Cal O_{\bP^n}(\l-\d+n+1)
\qaq C_{k,h}(\l)=H^0(\P,\cM_{k,h}(\l)).$$
From the exact sequence (5-1) we obtain the exact sequence
$$
0 \ra
\cM_{p,m+1}(\l) \ra
\cdots \ra
\cM_{p,1}(\l) \ra
\cM_{p,0}(\l) \ra 0,
\leqno(7-1)
$$
that induces the following complex
$$
0 \ra
C_{p,m+1}(\l) \ra
\cdots \ra
C_{p,1}(\l) \os{\phi}{\ra}
C_{p,0}(\l) \ra 0,
$$
Note that $\Coker \phi = B_p(\l)$ by Lem.(2-2) and (4-2).
We have the following commutative diagram:
$$
\matrix
\cdots & \ra &
C_{p,1}(\l)\ot \os{q+1}{\wedge} W & @>\phi>> &
C_{p,0}(\l)\ot \os{q+1}{\wedge} W & \ra &
B_p(\l)\ot \os{q+1}{\wedge} W & \ra & 0
\\
& & \downarrow & & \downarrow & & \downarrow &&
\\
\cdots & \ra &
C_{p+1,1}(\l)\ot \os{q}{\wedge} W & @>\phi>> &
C_{p+1,0}(\l)\ot \os{q}{\wedge} W & \ra &
B_{p+1}(\l)\ot \os{q}{\wedge} W & \ra & 0
\\
& & \downarrow & & \downarrow & & \downarrow &&
\\
\cdots & \ra &
C_{p+2,1}(\l)\ot \os{q-1}{\wedge} W & @>\phi>> &
C_{p+2,0}(\l)\ot \os{q-1}{\wedge} W & \ra &
B_{p+2}(\l)\ot \os{q-1}{\wedge} W & \ra & 0
\\
& & \downarrow & & \downarrow & & \downarrow &&
\\
& & \vdots & & \vdots & & \vdots &&
\endmatrix
$$
By an easy diagram chase, we see that the exactness of
$$ B_p(\l)\ot \os{q+1}{\wedge} W \to
B_{p+1}(\l)\ot \os{q}{\wedge} W \to
B_{p+2}(\l)\ot \os{q-1}{\wedge} W,$$
follows from the following.
\newline
\roster
\item
$
C_{p+2+a,a+1}(\l) \lra C_{p+2+a,a}(\l) \lra C_{p+2+a,a-1}(\l)
$
is exact for $1\leq\forall a \leq q-1$.
\item
$
C_{p+b,b}(\l)\ot \os{q+1-b}{\wedge}W \lra
C_{p+b+1,b}(\l)\ot \os{q-b}{\wedge}W \lra
C_{p+b+2,b}(\l)\ot \os{q-1-b}{\wedge}W
$
is exact for $\forall b \geq 0$.
\endroster
Note that (1) holds by trivial reason if $q=0$ or $q=1$.

\Lem 7-3. \it
Assume $s\geq 1$ and $e_1=\cdots=e_s$ and $p\geq 0$. (1) holds
in either of the following cases
\roster
\item"$(i)$"
$\ell\geq \d-n-1$ and $r+s\leq n+2$.
\item"$(ii)$"
$\ell\geq \d-n-1$ and $p\not= n-r-1$.
\endroster
\rm

\Lem 7-4. \it
Assume $s\geq 1$ and $p\geq 0$. (2) holds in either of the following cases
\roster
\item"$(i)$"
$q=0$ and $\md p+\ell\geq c$.
\item"$(ii)$"
$q=1$ and $\md p+\ell\geq 1+c$ and
$\md(p+1)+\ell\geq d_{max}+c$.
\item"$(iii)$"
$\md(r+p)+\ell\geq \d+q+c$ and $\ell\geq \d-n-1$.
\endroster
\rm \vskip 6pt

Before proving the lemmas, we finish the proof of Th.(7-1).
In case $(i)$, $(ii)$ and $(iii)$ it is a direct consequence of Lem.(7-4)
and (7-3). The case $(iv)$ is reduced to the case $(iii)$
by induction on $s$. For this we use the following commutative diagram
$$\matrix
0&&0&&0\\
\downarrow&&\downarrow&&\downarrow\\
\Rp {p-1} {\l+e_1} \ot \os{q+1}{\wedge} W &\to &
\Rp {p} {\l+e_1}\ot \os{q}{\wedge} W& \to&
\Rp {p+1} {\l+e_1}\ot \os{q-1}{\wedge} W\\
\downarrow&&\downarrow&&\downarrow\\
B_p(\l)\ot \os{q+1}{\wedge} W &\to &
B_{p+1}(\l)\ot \os{q}{\wedge} W& \to&
B_{p+2}(\l)\ot \os{q-1}{\wedge} W\\
\downarrow&&\downarrow&&\downarrow\\
\Rb p {\l}\ot \os{q+1}{\wedge} W &\to &
\Rb {p+1} {\l}\ot \os{q}{\wedge} W& \to&
\Rb {p+2} {\l}\ot \os{q-1}{\wedge} W\\
\downarrow&&\downarrow&&\downarrow\\
0&&0&&0\\
\endmatrix$$
where the notation is the same as in Lem.(5-3). The exactness of the
vertical sequences is a consequence of Lem.(5-3) and Th.(II)
(Here we use the additional assumption $\d+e_{max}-n-1> \ell$).
By the induction hypothesis we may assume the upper horizontal sequence
is exact. It remains to show the exactness of the lower horizontal sequence
$$\Rb p {\l}\ot \os{q+1}{\wedge} W \to
\Rb {p+1} {\l}\ot \os{q}{\wedge} W \to
\Rb {p+2} {\l}\ot \os{q-1}{\wedge} W$$
Letting $\Wb=\Im(W\to H^0(\Lb))$ and $I=\Ker(W\to H^0(\Lb))$, we have the
filtration $F^\cdot(\os{q}{\wedge} W)$ on $\os{q}{\wedge} W$ such that
$F^i/F^{i+1}\simeq (\os{q-i}{\wedge} \Wb)\ot (\os{i}{\wedge} I)$.
Since $I$ annihilates $\Rb p \l$, the above complex is filtered by
the above filtration and its graded quotients are the complexes
$$\Rb p {\l}\ot \os{q+1-i}{\wedge} \Wb\ot \os{i}{\wedge} I\to
\Rb {p+1} {\l}\ot \os{q-i}{\wedge} \Wb\ot\os{i}{\wedge} I \to
\Rb {p+2} {\l}\ot \os{q-i-1}{\wedge} \Wb\ot\os{i}{\wedge} I
\qfor 0\leq i\leq q.$$
These are exact by the induction hypothesis and this completes the proof.
\qed
\vskip 6pt

\demo{Proof of Lem.(7-3)}
The exact sequence (7-1) induces a spectral sequence
$$
E_1^{\alpha,\beta}=H^{\beta}(\cM_{k,m+1-\alpha}(\l))\Longrightarrow
H^{\alpha+\beta}=0.
$$
We want to show that $E^{\alpha,0}_2=0$ in case:
\roster
\item"$(*)$"
$p+3 \leq k \leq p+q+1$ and $k-(m+1-\alpha)=p+2$ ($\Longleftrightarrow$ $\alpha=p-k+ m+3$).
\endroster
Since $E^{\alpha,0}_{\infty}=0$, in order to show $E^{\alpha,0}_2=0$,
it suffices to show that
$$E^{\alpha-h-1,h}_1=H^h(\cM_{k,m+2+h-\alpha}(\l))=0
\qfor \forall h \geq 1.$$
In case $(*)$, putting $\ell'=\ell-\d+n+1$ we have
$$E_1^{\alpha-h-1,h}=H^h(\cM_{k,k-p+h-1}(\l))=
H^h(\os{m+2+p-(h+k)}{\wedge}\SS \ot \L^{p+r-h+1}
\ot \pi^*\Cal O_{\bP^n}(\l')).$$
We want to show that it vanishes assuming $k \geq p+3$ and $h \geq 1$.
We may suppose $h+k \leq m+2+p$ which implies $h\leq m-1$.
\vskip 4pt\noindent
$\underline{\text{Case}}$ $p+r-h+1 \leq -s$ \newline
We have $m-1\geq h\geq p+r+s+1>r+s$. The desired vanishing follows
from Lem(4-1)$(3)^*$.

\vskip 4pt\noindent
$\underline{\text{Case}}$ $p+r-h+1 \geq -s+1$ \newline
By the assumption $\l' \geq 0$ and by Lem(4-1)(1) and (3),
we have only to check the case $h=p+r+1=n$ and $\l'=0$ so that we are
concerned with the vanishing of $H^n(\os{n+s-k}{\wedge}\SS)$.
By Lem.(4-1)(4) this vanishies
if $n+s-k\leq n\Leftrightarrow s\leq k$. Since $k\geq p+3=n-r+2$, this holds
if $r+s\leq n+2$.
\vskip 5pt
This completes the proof of Lem.(7-3).
\qed
\enddemo

\demo{Proof of Lem.(7-4)}
By (4-2) we have
$$ C_{p,0}=H^0(\L^p\ot \pi^*\Cal O_{\bP^n}(\ell))=A_p(\ell) \qaq
 C_{p,1}=H^0(\L^{p-1}\ot\Sigma\ot \pi^*\Cal O_{\bP^n}(\ell)).
 \leqno(*)$$
Thus, Lem.(7-4) in case $(i)$ follows from Lem.(7-2).
In case $(ii)$ we need show the exactness of
$$C_{p,0}(\l)\ot \os{2}{\wedge}W \lra C_{p+1,0}(\l)\ot W \lra C_{p+2,0}(\l)$$
and the surjectivity of
$C_{p+1,1}(\l)\ot W \lra C_{p+2,1}(\l)$.
By $(*)$ the first assertion follows from Lem.(7-2). To show the
second assertion we recall the exact sequences
$$
0 \ra
\Cal O_{\P} \ra
\Sigma_{\L}(-\log \P_*) \ra
T_{\P}(-\log\P_*) \ra 0,
$$
$$
0 \ra
T_{\P/\P^n}(-\log \P_*) \ra
T_{\P}(-\log \P_*) \ra
\pi^*T_{\P^n}\to 0,
$$
$$
0\ra \Cal O_{\P} \ra
\pi^*\Cal E_0^*\ot \Cal L \op \Cal O_{\Bbb P}^{\oplus s}
\ra T_{\Bbb P/\Bbb P^n}(-\log \P_*)\ra 0,
$$
$$
0\ra \Cal O_{\P^n} \ra
\Cal O_{\Bbb P^n}(1)^{\oplus n+1}
\ra T_{\Bbb P^n}\ra 0.
$$
Noting $H^1(\L^\nu\ot\pi^*\Cal O_{\bP^n}(\ell))=0$ for $\forall \nu\geq 0$,
the assertion follows from the surjectivity of
$$ H^0(\L^p\ot\pi^*\Cal O_{\bP^n}(\ell))\ot W\to
H^0(\L^{p+1}\ot\pi^*\Cal O_{\bP^n}(\ell)),$$
$$ H^0(\L^{p+1}\ot\pi^*\Cal E_0^*(\ell))\ot W \to
 H^0(\L^{p+2}\ot\pi^*\Cal E_0^*(\ell))$$
which is a consequence of Lem.(7-2) under the assumption of $(ii)$.
\par

Finally we show Lem.(7-4) in case $(iii)$.
We denote $\Om_{\P}^1(\log \P_*)$ by $\Om$ simply.
We want to show that
$$
\multline
H^0(\os{m+1-b}{\wedge}\SS \ot \L^{r+p}\ot \pi^*\Cal O_{\bP^n}(\l')) \ot \os{q+1-b}{\wedge}W \\ \ra
H^0(\os{m+1-b}{\wedge}\SS \ot \L^{r+p+1}\ot \pi^*\Cal O_{\bP^n}(\l')) \ot \os{q-b}{\wedge}W \\ \ra
H^0(\os{m+1-b}{\wedge}\SS \ot \L^{r+p+2}\ot \pi^*\Cal O_{\bP^n}(\l')) \ot \os{q-1-b}{\wedge}W
\endmultline
$$
is exact for $0 \leq \forall b \leq q$, where $\ell'=\ell-\d+n+1$.
If $\sharp >0$ and $\ell'\geq 0$, the following exact sequence
$$
0 \ra
\Om^{\cdot}\ot \L^{\sharp}\ot \pi^*\Cal O_{\bP^n}(\l') \ra
\os{\cdot}{\wedge}\SS \ot \L^{\sharp}\ot \pi^*\Cal O_{\bP^n}(\l') \ra
\Om^{\cdot-1}\ot \L^{\sharp}\ot \pi^*\Cal O_{\bP^n}(\l') \ra 0
$$
remains exact after taking $H^0(\quad)$ since
$H^1(\L^{\sharp}\ot\Om^{\cdot}\ot \pi^*\Cal O_{\bP^n}(\l'))=0$, which we can see from
the proof of Lem.(4-1) (cf. Claim below (4-7)).
Thus it suffices to show that the following sequence is exact for all
$t$, $b$ such that $m-b \leq t \leq m-b+1$ and $0 \leq b \leq q$:
$$
\multline
H^0(\Om^t \ot \L^{r+p}\ot \pi^*\Cal O_{\bP^n}(\l')) \ot \os{q+1-b}{\wedge}W \\ \ra
H^0(\Om^t \ot \L^{r+p+1}\ot \pi^*\Cal O_{\bP^n}(\l')) \ot \os{q-b}{\wedge}W \\ \ra
H^0(\Om^t \ot \L^{r+p+2}\ot \pi^*\Cal O_{\bP^n}(\l')) \ot \os{q-1-b}{\wedge}W.
\endmultline
$$
By (4-5), there is a filtration $F^{\cdot}$ of
$\Om^t$ such that
$$\text{Gr}^u_F(\Om^t)
=\pi^*\Om^u_{\P^n}\ot (\os{t-u}{\us{i=0}{\op}}[\os{i}{\wedge}\pi^*\E_0
\ot \L^{-i}]^{\binom{s-1}{t-u-i}}),$$
where $(u,i)$ runs over
$0 \leq u \leq n$ and $0\leq i \leq \min\{t-u,r\}$.
Since $H^1(\L^{r+\sharp}\ot \text{Gr}^{\cdot}_F\Om^t\ot \pi^*\Cal O_{\bP^n}(\l'))=0$
for $\sharp \geq 0$ and $\ell'\geq 0$, it suffices to show that
$$
\multline
H^0( \L^{r+p-i}\ot \pi^*(\Om^u_{\P^n}\ot \os{i}{\wedge}\E_0(\l')) \ot \os{q+1-b}{\wedge}W \\ \ra
H^0( \L^{r+p-i+1}\ot \pi^*(\Om^u_{\P^n}\ot \os{i}{\wedge}\E_0(\l')) \ot \os{q-b}{\wedge}W \\ \ra
H^0( \L^{r+p-i+2}\ot \pi^*(\Om^u_{\P^n}\ot \os{i}{\wedge}\E_0(\l')) \ot \os{q-1-b}{\wedge}W.
\endmultline
$$
$$
\text{ is exact for $\forall b,u,i$ such that}
 \cases
 0 \leq b \leq q \\
 0 \leq u \leq n \\
 0\leq i\leq \min\{t-u,r\}\\
 \endcases
$$
Finally, by the exact sequence
$$
0 \ra \Cal O_{\Bbb P^n}(-n-1)
  \ra \Cal O_{\Bbb P^n}(-n)^{\op n+1}
  \ra \cdots
  \ra \Cal O_{\Bbb P^n}(-u-1)^{\op{n+1 \choose n-u}}
  \ra \Omega^u_{\Bbb P^n}
  \ra 0,
$$
we can reduce the assertion to show that
$$
\multline
H^0( \L^{r+p-i+j}\ot \pi^*(\os{i}{\wedge}\E_0(\l'-u-j-1)) \ot \os{q+1-b-j}{\wedge}W \\ \ra
H^0( \L^{r+p-i+j+1}\ot \pi^*(\os{i}{\wedge}\E_0(\l'-u-j-1)) \ot \os{q-b-j}{\wedge}W \\ \ra
H^0( \L^{r+p-i+j+2}\ot \pi^*(\os{i}{\wedge}\E_0(\l'-u-j-1)) \ot \os{q-1-b-j}{\wedge}W
\endmultline
$$
$$
\text{ is exact for $\forall b,u,i,j$ such that} \quad
\cases
 0 \leq b \leq q\\
 0 \leq u \leq n        \\
 0\leq i\leq \min\{t-u,r\}\\
 0 \leq j \leq n-u \\
\endcases
$$
Let $\delta_i$ be the minimal degree of line bundles which are direct summands
of $\os{i}{\wedge}\E_0$. Then by Lem.(7-2), the above holds if $p \geq 0$ and
$$
\md(r+p-i+j)+(\l'+\delta_i-u-j-1) \geq q-b-j+c
$$
for $\forall b,u,i,j$ as above. By noting $\delta_i\geq \md i$, it is easy to
see that this holds under the assumption
$\md(r+p)+\l-\d\geq q+c$.
This completes the proof.
\qed
\enddemo


\vskip 20pt

\input amstex
\documentstyle{amsppt}
\hsize=16cm
\vsize=23cm

\def\Th#1.{\vskip 6pt \medbreak\noindent{\bf Theorem(#1).}}
\def\Cor#1.{\vskip 6pt \medbreak\noindent{\bf Cororally(#1).}}
\def\Conj#1.{\vskip 6pt \medbreak\noindent{\bf Conjecture(#1).}}
\def\Pr#1.{\vskip 6pt \medbreak\noindent{\bf Proposition(#1).}}
\def\Lem#1.{\vskip 6pt \medbreak\noindent{\bf Lemma(#1).}}
\def\Rem#1.{\vskip 6pt \medbreak\noindent{\it Remark(#1).}}
\def\Fact#1.{\vskip 6pt \medbreak\noindent{\it Fact(#1).}}
\def\Claim#1.{\vskip 6pt \medbreak\noindent{\it Claim(#1).}}
\def\Def#1.{\vskip 6pt \medbreak\noindent{\bf Definition\bf(#1)\rm.}}

\def\qwith{\quad\hbox{with }}
\def\mathrm#1{\rm#1}

\def\isom{@>\cong>>}

\def\Ext{{\text{\rm{Ext}}}}

\def\dim{{\operatorname{dim}}}

\def\Coker{{\text{\rm Coker}}}
\def\dim{\hbox{\rm dim}}
\def\det{\hbox{\rm det}}

\def\Im{\hbox{\rm Im}}
\def\Ker{\hbox{\rm Ker}}
\def\Coker{\hbox{\rm Coker}}
\def\min{\hbox{\rm min}}

\def\Hom{\hbox{\rm{Hom}}}

\def\sign{\hbox{\mathrm{sign}}}
\def\Res#1{\hbox{\mathrm{Res}}_{#1}}

\def\P{{\Bbb{P}}}
\def\bP{{\Bbb{P}}}

\def\Q{{\Bbb{Q}}}

\def\cHom{{\Cal{H}}om}
\def\cL{{\Cal{L}}}
\def\L{{\Cal{L}}}
\def\cE{{\Cal{E}}}
\def\E{{\Cal{E}}}
\def\cF{{\Cal{F}}}
\def\cO{{\Cal{O}}}
\def\cX{{\Cal{X}}}

\def\cZ{{\Cal{Z}}}

\def\cM{{\Cal{M}}}

\def\sE#1#2#3{E_{#1}^{#2,#3}}

\def\l{\ell}
\def\d{{\bold{d}}}
\def\e{{\bold{e}}}

\def\lra{\longrightarrow}
\def\ra{\rightarrow}

\def\hra{\hookrightarrow}

\def\ot{\otimes}
\def\op{\oplus}

\def\us#1#2{\underset{#1}\to{#2}}
\def\os#1#2{\overset{#1}\to{#2}}

\def\PE{\P(\cE)}
\def\qaq{\quad\hbox{ and }\quad}
\def\qfor{\quad\hbox{ for }}
\def\qif{\quad\hbox{ if }}
\def\ld{\lambda}
\def\spa{\hbox{ }}
\def\scs{\spa : \spa}
\def\SS{\Sigma^*}
\def\sp{\hbox{}}

\def\Sp{\Sigma'}
\def\Sb{\overline{\Sigma}}
\def\SSp{\Sigma^{'*}}

\def\taub{\overline{\tau}}
\def\Lb{\overline{\Cal L}}
\def\Pb{\overline{\Bbb P}}
\def\Rp#1#2{B'_{#1}(#2)}
\def\Rb#1#2{\overline{B}_{#1}(#2)}
\def\eb{\overline{\bold e}}
\def\pg#1{{\frak S}_{#1}}
\def\ul#1{\underline{#1}}

\def\ba{\bold{a}}
\def\bb{\bold{b}}
\def\ua{\underline{a}}
\def\ub{\underline{b}}

\def\ue{\underline{e}}
\def\ud{\underline{d}}

\def\md{\delta_{\text{min}}}
\def\Om{\Omega}
\def\ld{\lambda}
\def\Wb{\overline{W}}

\def\FS#1{F_S^{#1}}
\def\WS#1{\Omega_S^{#1}}

\def\otO{\otimes_{\cO}}
\def\cU{{\Cal U}}

\def\onab{\overline{\nabla}}

\def\HcU#1#2{H^{#1,#2}(\cU/S)}
\def\HccU#1#2{H_{\cO,c}^{#1,#2}(\cU/S)}

\def\HcU#1#2{H^{#1,#2}(\cU/S)}

\def\HX#1#2{H^{#1}(X,{#2})}

\def\Xx{X_x}
\def\Zx{Z_x}
\def\Ux{U_x}

\def\Zxx{Z_{x}}
\def\Zst{Z}
\def\cZst{\cZ}
\def\Zast{Z^{(\alpha)}}
\def\Wst{W}
\def\Yst{Y}

\def\rocXZ{\kappa_{(\cX,\cZ)}}

\def\roolog{\kappa_0^{log}}

\def\roxlog{\kappa_x^{log}}

\def\TS{\Theta_S}
\def\TxS{T_x S}
\def\ToS{T_0 S}

\def\WS#1{\Omega_S^{#1}}
\def\WX#1{\Omega_X^{#1}}

\def\TX{T_X}

\def\TXZx{T_{\Xx}(-\log \Zxx)}

\def\WXZ#1{\Omega_X^{#1}(\log \Zst)}
\def\WXZa#1{\Omega_X^{#1}(\log \Zast)}
\def\WcXZ#1{\Omega_{\cX/S}^{#1}(\log \cZst)}

\def\WcXkZ#1{\Omega_{\cX/k}^{#1}(\log \cZst)}
\def\WcXCZ#1{\Omega_{\cX/\Bbb C}^{#1}(\log \cZst)}
\def\WPnY#1{\Omega_{\Bbb P^n}^{#1}(\log \Yst)}
\def\TXZ{T_X(-\log \Zst)}
\def\TcXZS{T_{\cX/S}(-\log \cZst)}

\def\fXZ#1{\phi_{X,Z}^{#1}}
\def\cXZ{\psi_{(X,Z)}}

\def\cXZx{\psi_{(\Xx,\Zx)}}

\def\eXZ{\eta_{(X,Z)}}

\def\ccXZS{c_S(\cX,\cZ)}

\def\dpXZ#1{d\rho^{#1}_{X,Z}}
\def\scs{\hbox{ }:\hbox{ }}
\def\onab{\overline{\nabla}}
\def\tonab{@>\onab>>}

\def\spa{\hbox{ }}

\def\otk{\otimes_k}
\def\Pol{P}
\def\h#1#2{h_{#1}(#2)}

\def\Szar{S_{zar}}

\head \S8. Infinitesimal Torelli for open complete intersections \endhead
\vskip 8pt

Let the notation and the assumption be as in Def.(1-2).
The main result in this section is the infinitesimal Torelli
for the pair $(X,Z)$, which concerns the injectivity of the following map
$$\dpXZ q \scs \HX 1 {\TXZ} \to
\Hom(\HX {n-r-q} {\WXZ q},\HX {n-r-q+1} {\WXZ {q-1}}) $$
where $1\leq q\leq n-r$ and $\WXZ q$ is the sheaf of algebraic differential
$q$-forms on $X$ with logarithmic poles along $Z$ and $\TXZ$ is the
$\cO_X$-dual of $\WXZ 1$.
The above map is induced by the cup product and the contraction
$\WXZ q \ot \TXZ \to \WXZ {q-1}$.

\Th 8-1. \it Assume $\md(n-r-q)+\d+\e \geq n-1$ and
$\md(q-1)+\d\geq n-1$. Then $\dpXZ q$ is injective.
\rm\demo{Proof}
By noting $\WX {n-r}\simeq \cO_X(\d-n-1)$, the dual of $\dpXZ q$ is
identified with the map
$$ \HX {n-r-q} {\WXZ q}\otimes {\HX {n-r-q+1} {\WXZ {q-1}}}^*\to
H^{n-r-1}(X,\WXZ 1\otimes\cO_X(\d-n-1)).$$
By Th.(I) and Th.(II) in \S1 it is identified with the multiplication of
Jacobian rings
$$ B_{n-r-q}(\d+\e-n-1)\otimes B_{q-1}(\d-n-1) \to B_{n-r-1}(2(\d-n-1)+\e).$$
By definition the condition of Th.(8-1) implies that every (bi)homogeneous
polynomial appearing in the Jacobian rings on the left hand side has a
non-negative degree. Hence the above map is surjective under the assumption.
\qed
\enddemo


\vskip 20pt

\input amstex
\documentstyle{amsppt}
\hsize=16cm
\vsize=23cm

\def\Th#1.{\vskip 6pt \medbreak\noindent{\bf Theorem(#1).}}
\def\Cor#1.{\vskip 6pt \medbreak\noindent{\bf Cororally(#1).}}
\def\Conj#1.{\vskip 6pt \medbreak\noindent{\bf Conjecture(#1).}}
\def\Pr#1.{\vskip 6pt \medbreak\noindent{\bf Proposition(#1).}}
\def\Lem#1.{\vskip 6pt \medbreak\noindent{\bf Lemma(#1).}}
\def\Rem#1.{\vskip 6pt \medbreak\noindent{\it Remark(#1).}}
\def\Fact#1.{\vskip 6pt \medbreak\noindent{\it Fact(#1).}}
\def\Claim#1.{\vskip 6pt \medbreak\noindent{\it Claim(#1).}}
\def\Def#1.{\vskip 6pt \medbreak\noindent{\bf Definition\bf(#1)\rm.}}

\def\qwith{\quad\hbox{with }}
\def\mathrm#1{\rm#1}

\def\isom{@>\cong>>}

\def\Ext{{\text{\rm{Ext}}}}

\def\dim{{\operatorname{dim}}}

\def\Coker{{\text{\rm Coker}}}
\def\dim{\hbox{\rm dim}}
\def\det{\hbox{\rm det}}

\def\Im{\hbox{\rm Im}}
\def\Ker{\hbox{\rm Ker}}
\def\Coker{\hbox{\rm Coker}}
\def\min{\hbox{\rm min}}

\def\Hom{\hbox{\rm{Hom}}}

\def\sign{\hbox{\mathrm{sign}}}
\def\Res#1{\hbox{\mathrm{Res}}_{#1}}

\def\P{{\Bbb{P}}}
\def\bP{{\Bbb{P}}}

\def\Q{{\Bbb{Q}}}

\def\cHom{{\Cal{H}}om}
\def\cL{{\Cal{L}}}
\def\L{{\Cal{L}}}
\def\cE{{\Cal{E}}}
\def\E{{\Cal{E}}}
\def\cF{{\Cal{F}}}
\def\cO{{\Cal{O}}}
\def\cX{{\Cal{X}}}

\def\cZ{{\Cal{Z}}}

\def\cM{{\Cal{M}}}

\def\sE#1#2#3{E_{#1}^{#2,#3}}

\def\l{\ell}
\def\d{{\bold{d}}}
\def\e{{\bold{e}}}

\def\lra{\longrightarrow}
\def\ra{\rightarrow}

\def\hra{\hookrightarrow}

\def\ot{\otimes}
\def\op{\oplus}

\def\us#1#2{\underset{#1}\to{#2}}
\def\os#1#2{\overset{#1}\to{#2}}

\def\PE{\P(\cE)}
\def\qaq{\quad\hbox{ and }\quad}
\def\qfor{\quad\hbox{ for }}
\def\qif{\quad\hbox{ if }}
\def\ld{\lambda}
\def\spa{\hbox{ }}
\def\scs{\spa : \spa}
\def\SS{\Sigma^*}
\def\sp{\hbox{}}

\def\Sp{\Sigma'}
\def\Sb{\overline{\Sigma}}
\def\SSp{\Sigma^{'*}}

\def\taub{\overline{\tau}}
\def\Lb{\overline{\Cal L}}
\def\Pb{\overline{\Bbb P}}
\def\Rp#1#2{B'_{#1}(#2)}
\def\Rb#1#2{\overline{B}_{#1}(#2)}
\def\eb{\overline{\bold e}}
\def\pg#1{{\frak S}_{#1}}
\def\ul#1{\underline{#1}}

\def\ba{\bold{a}}
\def\bb{\bold{b}}
\def\ua{\underline{a}}
\def\ub{\underline{b}}

\def\ue{\underline{e}}
\def\ud{\underline{d}}

\def\md{\delta_{\text{min}}}
\def\Om{\Omega}
\def\ld{\lambda}
\def\Wb{\overline{W}}

\def\FS#1{F_S^{#1}}
\def\WS#1{\Omega_S^{#1}}

\def\otO{\otimes_{\cO}}
\def\cU{{\Cal U}}

\def\onab{\overline{\nabla}}

\def\HcU#1#2{H^{#1,#2}(\cU/S)}
\def\HccU#1#2{H_{\cO,c}^{#1,#2}(\cU/S)}

\def\HcU#1#2{H^{#1,#2}(\cU/S)}

\def\HX#1#2{H^{#1}(X,{#2})}

\def\Xx{X_x}
\def\Zx{Z_x}
\def\Ux{U_x}

\def\Zxx{Z_{x}}
\def\Zst{Z}
\def\cZst{\cZ}
\def\Zast{Z^{(\alpha)}}
\def\Wst{W}
\def\Yst{Y}

\def\rocXZ{\kappa_{(\cX,\cZ)}}

\def\roolog{\kappa_0^{log}}

\def\roxlog{\kappa_x^{log}}

\def\TS{\Theta_S}
\def\TxS{T_x S}
\def\ToS{T_0 S}

\def\WS#1{\Omega_S^{#1}}
\def\WX#1{\Omega_X^{#1}}

\def\TX{T_X}

\def\TXZx{T_{\Xx}(-\log \Zxx)}

\def\WXZ#1{\Omega_X^{#1}(\log \Zst)}
\def\WXZa#1{\Omega_X^{#1}(\log \Zast)}
\def\WcXZ#1{\Omega_{\cX/S}^{#1}(\log \cZst)}

\def\WcXkZ#1{\Omega_{\cX/k}^{#1}(\log \cZst)}
\def\WcXCZ#1{\Omega_{\cX/\Bbb C}^{#1}(\log \cZst)}
\def\WPnY#1{\Omega_{\Bbb P^n}^{#1}(\log \Yst)}
\def\TXZ{T_X(-\log \Zst)}
\def\TcXZS{T_{\cX/S}(-\log \cZst)}

\def\fXZ#1{\phi_{X,Z}^{#1}}
\def\cXZ{\psi_{(X,Z)}}

\def\cXZx{\psi_{(\Xx,\Zx)}}

\def\eXZ{\eta_{(X,Z)}}

\def\ccXZS{c_S(\cX,\cZ)}

\def\dpXZ#1{d\rho^{#1}_{X,Z}}
\def\scs{\hbox{ }:\hbox{ }}
\def\onab{\overline{\nabla}}
\def\tonab{@>\onab>>}

\def\spa{\hbox{ }}

\def\otk{\otimes_k}
\def\Pol{P}
\def\h#1#2{h_{#1}(#2)}

\def\Szar{S_{zar}}

\head \S9. Explicit bound for Nori's connectivity \endhead
\vskip 8pt

In this section we deduce Th.(0-1) from Th.(9-1) and Th.(9-3),
the symmetrizer lemmas for open complete intersections.
Let the assumption be as in \S1. We fix a non-singular algebraic
variety $S$ over $k$ and the following schemes over $S$
$$\Bbb P^n_S \hookleftarrow \cX \overset{i}\to\hookleftarrow
\cZ=\underset{1\leq j\leq s}\to{\cup} \cZ_j\leqno(9-1)$$
whose fibers are as in Def.(1-2).
Let $f:\cX\to S$ be the natural morphism and write $\cU=\cX \setminus \cZ$.
For integers $p,q$ we introduce the following sheaf on $\Szar$
$$\HcU p q=R^q f_* \WcXZ p,$$
where $\WcXZ p=\os{p}{\wedge}\WcXZ 1$ with $\WcXZ 1$, the sheaf of relative
differentials on $\cX$ over $S$ with logarithmic poles along $\cZ$.
We assume $s\geq 1$. Then the Lefschetz theory implies $\HcU p q=0$
if $p+q\not=n-r$. Weconsider the following Koszul complex
$$ \WS {q-1}\otO \HcU {a+2}{b-2} @>\onab>>
 \WS {q}\otO \HcU {a+1}{b-1} @>\onab>>
 \WS {q+1}\otO \HcU {a}{b}.\leqno(9-2) $$
Here $\onab$ is induced by the Kodaira-Spencer map
$$ \rocXZ\scs \TS \to R^1f_* \TcXZS,$$
with $\TS=\cHom_{\cO_S}(\WS 1,\cO_S)$ and
$\TcXZS=\cHom_{\cO_{\cX}}(\WcXZ 1,\cO_{\cX})$, and the map
$$ R^1f_* \TcXZS \otimes R^{b-1} f_* \WcXZ {a+1} \to R^b f_* \WcXZ a$$
induced by the cup product and
$\TcXZS \otimes \WcXZ {a+1} \to \WcXZ a$, the contraction.

\Th 9-1. \it
Let $c=\ccXZS$ be as in Def.(9-1) below and assume $n-r\geq 2$.
Assume also that $(*)$ either $a<n-r-1$ or $r+s\leq n$.
Then the complex (9-2) with $a+b=n-r$ is exact under one of the following
conditions
\roster
\item"$(i)$"
$a\geq 0$, $q=0$ and $\md a+\d\geq c+n+1$.
\item"$(ii)$"
$a\geq 0$, $q=1$ and $\md a+\d\geq c+n+2$ and
$\md(a+1)+\d \geq c+n+1+d_{max}$.
\item"$(iii)$"
$a\geq 0$, $\md(r+a)\geq q+c+n+1$ and $r+s\leq n+2$.
\item"$(iv)$"
$a\geq 0$, $\md(r+a)\geq q+c+n+1$ and $a<n-r-\frac{q}{2}$.
\endroster

\Def 9-1. \it For $x\in S$ let $\Ux\subset\Xx\supset\Zx$ denote the fibers of
the family (9-1) and let
$$ \roxlog \scs \TxS \to H^1(\Xx,\TXZx)
\quad (resp. \spa \cXZx: B_1(0) \to H^1(\Xx,\TXZx))$$
be the Kodaira-Spencer map (resp. the map in Th.(I)(2) for $(\Xx,\Zx)$).
We define
$$\ccXZS=\underset{\spa x\in S}\to{\max}
\{\dim_{k}(\Im(\cXZx)/\Im(\cXZx)\cap \Im(\roxlog))\}.$$
If $n-r\geq 2$ and $\Xx$ is not a $K3$ surface,
$\cXZx$ is surjective so that
$$\ccXZS=\underset{x\in S}\to{\max}\{\dim_{k(x)}
(\Coker(\rocXZ)\otimes_{\cO_S} k(x))\} .$$
\vskip 6pt\rm

Now we prove Th.(9-1).
We fix $0\in S$ and let $X\supset Z$ be the fiber over $0$ of $\cX\supset\cZ$.
By Th.(I) and (II) the assumption $(*)$ implies that the dual of the fiber over
$0$ of the complex (9-2) is identified with
$$
B_a(\l_0) \ot \os{q+1}{\wedge} \ToS \lra
  B_{a+1}(\l_0) \ot \os{q}{\wedge} \ToS\lra
B_{a+2}(\l_0) \ot \os{q-1}{\wedge} \ToS \quad
(\l_0:=\d-n-1)
\leqno(**)
$$
where $B_*(\ell)$ denotes the Jacobian ring for $(X,Z)$
and the maps are induced by the composite map
$$ \rho\scs \ToS @>{\roolog}>> H^1(X,\TXZ)_{alg}
@>{\cXZ^{-1}}>> B_1(0)\quad \text{(cf. Th.(I))}$$
and the multiplication on the Jacobian rings.
Let $V=\Im(\rho)\subset B_1(0)$ and $K=\Ker(\rho)$. By definition $V$ is of
codimension$\leq c$ in $B_1(0)$. We have the filtration
$F^\nu (\os{q}{\wedge} \ToS)\subset \os{q}{\wedge} \ToS$ such that
$$F^\nu (\os{q}{\wedge} \ToS)/F^{\nu-1} (\os{q}{\wedge} \ToS) \isom
\os{\nu}{\wedge}K\otimes \os{q-\nu}{\wedge} V.$$
The filtration induces a filtration of the complex $(**)$ so that
it suffices to show the exactness of $(**)$ with $\ToS$ replaced by $V$.
Then it follows from Th.(III).
\qed
\vskip 6pt

\Th 9-2. \it Let $c=\ccXZS$ be as in Def.(9-1). Assume $n-r\geq 2$
and that $S$ is affine.
Then
$$H^b(\cX,\WcXkZ a)=0
\qif
s\leq n-r+2,\spa b\leq n-r-1,\spa
\md(n-1-b)\geq a+b+1+r+c,$$
where $\WcXkZ \cdot$ is the sheaf of differential forms of $\cX$ over $k$
with logarithmic poles along $\cZ$.
\rm\vskip 6pt
\demo{Proof}
Filter $\WcXkZ a$ by the subsheaves
$$ \FS q \WcXkZ a= \Im(f^*\WS q\ot\WcXkZ {a-q} \to\WcXkZ a)$$
so that
$$ Gr_{F_S}^q\WcXkZ a =f^*\WS q\ot\WcXZ {a-q}.$$
The filtration gives rise to the spectral sequence
$$ E_1^{q,p}=H^{q+p}(Gr_{F_S}^q\WcXkZ a)=\WS q\ot \HcU {a-q}{q+p}
\Rightarrow H^{q+p}(\cX,\WcXkZ a).$$
By the Lefschetz theory $E_1^{q,b-q}=0$ unless $a+b-q=n-r$ in which case
$E_2^{q,b-q}$ is computed as the cohomology of the complex (9-2).
Th.(9-1) implies that $E_2^{q,b-q}=0$ if $s\leq n-r+2$ and $a-q-1\geq 0$ and
$\md(r+a-q-1)\geq q+c+n+1$, which is in case $a+b-q=n-r$ equivalent to
the assumption of Th.(9-2).
\qed
\enddemo
\vskip 5pt

Now the first vanishing of Th.(0-1) is an easy consequence of Th.(9-2)
since the vanishing of $F^{\ell} H^t(\cU,\Bbb C)$ is
reduced to that of $\Bbb H^t(\cX,\WcXCZ {\geq \l})$
by [D1, Pr.3.1.8].
\vskip 6pt

In order to show the second vanishing of Th.(0-1) we consider
the family (9-1) assuming $s=1$. For integers $a,b$ we write
$\HccU a b =R^b f_* \Omega_{(\cX,\cZ)/S}^a$
where $\Omega^a_{(\cX,\cZ)/S}$ is defined by the exact sequence
$$0\to\Omega^a_{(\cX,\cZ)/S}\to\Omega^a_{\cX/S}\to i_*\Omega^a_{\cZ/S}\to 0.$$
By the Lefschetz theory and the Serre duality
$$ \HccU a b=0 \qif a+b\not=n-r \qaq
\HccU a b = {\HcU b a}^*.$$
We consider the complex
$$\WS {q-1}\otO \HccU {a+2} {b-2} \tonab
\WS {q}\otO \HccU {a+1} {b-1} \tonab
\WS {q+1}\otO \HccU {a} {b}\leqno(9-3)$$
where the maps are induced by the Gauss-Manin connection.
By the same argument as the proof of Th.(9-1) we can show the following.

\Th 9-3. \it Assume $n-r\geq 2$ and $s=1$ and write $e=e_1$.
Let $c=\ccXZS$ be as in Def.(9-1).
The complex (9-3) with $a+b=n-r$ is exact under one of the following
conditions
\roster
\item"$(i)$"
$a\geq 0$, $q=0$ and $\md a+\d+e \geq c+n+1$.
\item"$(ii)$"
$a\geq 0$, $q=1$ and $\md a+\d+e\geq c+n+2$ and
$\md(a+1)+\d+e \geq c+n+1+d_{max}$.
\item"$(iii)$"
$a\geq 0$, $\md(r+a)+e\geq q+c+n+1$.
\endroster
\rm\vskip 6pt

\Th 9-4. \it Let $c=\ccXZS$ be as in Def.(9-1). Assume $n-r\geq 2$
and that $S$ is affine. Then
$$H^b(\cX,\Omega^a_{(\cX,\cZ)/k})=0
\qif
b\leq n-r-1
\qaq
\md(n-1-b)+e\geq a+b+1+r+c$$
where $\Omega^a_{(\cX,\cZ)/k}$ is defined by the exact sequence
$0\to\Omega^a_{(\cX,\cZ)/k}\to\Omega^a_{\cX/k}\to i_*\Omega^a_{\cZ/k}\to 0.$
\rm\vskip 6pt
\demo{Proof}
Filter $\Omega^a_{(\cX,\cZ)/k}$ by the subsheaves
$$ \FS q \Omega^a_{(\cX,\cZ)/k}=
\Im(f^*\Omega_S^q\otimes\Omega^{a-q}_{(\cX,\cZ)/k}\to\Omega^a_{(\cX,\cZ)/k})$$
so that
$$ Gr_{F_S}^q\Omega^a_{(\cX,\cZ)/k}=
f^*\Omega_S^q\otimes\Omega^{a-q}_{(\cX,\cZ)/S}.$$
The rest of the argument is the same as the proof of Th.(9-2).
\qed
\enddemo
\vskip 8pt

As is shown in [N], the vanishing of $F^{\ell} H^t(\cX,\cZ,\Bbb C)$ is
reduced to that of $\Bbb H^t(\cX,\Omega^{\geq \ell}_{(\cX,\cZ)/\Bbb C})$.
Thus the second vanishing of Th.(0-1) is an easy consequence of Th.(9-4).


\vskip 20pt

\input amstex
\documentstyle{amsppt}
\hsize=16cm
\vsize=23cm

\document
\NoBlackBoxes
\nologo


\enddocument